\documentclass[preprint,review,3p,times]{elsarticle}

\usepackage{amssymb}
\usepackage{amsmath}

\usepackage{lineno}

\usepackage{float}
\usepackage{subcaption}
\usepackage{placeins}

\usepackage{tikz}
\usepackage{tikz-3dplot}
\usetikzlibrary{calc,fadings,decorations.pathreplacing,patterns,patterns.meta}

\usetikzlibrary{shapes.geometric, arrows}

\tikzstyle{io} = [trapezium,
trapezium left angle=70, 
trapezium right angle=110, 
minimum width=2cm, 
minimum height=1cm,
text centered, 
text width=2cm,
draw=black, 
fill=blue!30]

\tikzstyle{process} = [rectangle, 
minimum width=3.5cm, 
minimum height=1cm, 
text centered, 
text width=3.5cm, 
draw=black, 
fill=orange!30]

\tikzstyle{arrow} = [thick,->,>=stealth]

\tikzset{%
  >=latex, 
  inner sep=0pt,%
  outer sep=2pt,%
  mark coordinate/.style={inner sep=0pt,outer sep=0pt,minimum size=3pt,
    fill=black,circle}%
}

\journal{Journal of Sound and Vibration}

\begin{document}
\begin{frontmatter}

\title{A space-time Galerkin boundary element method for aeroacoustic scattering}

\author[label1]{Maks Groom\corref{cor1}}
\author[label1]{Beckett Zhou}

\affiliation[label1]{organization={Daniel Guggenheim School of Aerospace Engineering, Georgia Institute of Technology},
city={Atlanta},
state={Georgia},
country={USA}}

\cortext[cor1]{Corresponding author.\\ \noindent \hspace*{1.5em} \textit{E-mail address:} mgroom@gatech.edu (M. Groom).}

\begin{abstract}
Acoustic scattering by vehicle surfaces can have significant effects on overall noise levels. In this paper, we present a space-time Galerkin time-domain boundary element method (TDBEM) that offers several distinct advantages over contemporary scattering methods for prediction of acoustic scattering and shielding of complex aeroacoustic sources such as propellers and rotors. The time-domain approach allows efficient simulation of transient, rotating, and broadband noise sources, while the Galerkin formulation is robust and unconditionally stable without any tuned numerical parameters. The main challenge of the Galerkin approach, namely the numerically difficult double space-time integration, is resolved through an efficient decomposition-based quadrature procedure. We present three cases with analytical solutions to validate the method and study its numerical properties, demonstrating excellent agreement for scattering and shielding by a variety of different geometries. We then apply the TDBEM to a trailing edge-mounted propeller case, comparing the numerical predictions with experimental measurements. The results demonstrate good agreement between predicted and measured scattering and shielding in a practical application case.
\end{abstract}

\begin{keyword}
Aeroacoustic scattering \sep time-domain boundary element method \sep propeller noise



\end{keyword}

\end{frontmatter}

\section{Introduction}\label{sec1}
Vehicle noise is an important design consideration for many aircraft to mitigate adverse health and environmental impacts and improve public acceptability~\cite{rizzi2020urban}. Scattering and shielding by vehicle surfaces can significantly modify the radiated noise from aeroacoustic sources, affecting the overall aircraft noise levels. These effects are not typically captured in widely-used integral-based acoustic analogy methods~\cite{dunn2019acoustic}. In some cases, scattering can be directly captured with acoustic analogies if high-fidelity compressible computational fluid dynamics (CFD) simulations are used to obtain an acoustically resolved flow solution, and all noise sources are subsequently included in the integration~\cite{romani2017application}. However, resolution of acoustic fluctuations within conventional CFD simulations is severely limited by numerical artifacts and computational constraints~\cite{tam2012computational}. Additionally, many formulations that are inherently unable to resolve acoustic fluctuations, including incompressible CFD~\cite{yangzhou2023aeroacoustic}, are frequently applied to noise prediction in combination with acoustic analogies. It is therefore important to develop computationally efficient scattering prediction tools to enable accurate assessment of vehicle noise early in the design process.

There are several approaches for including scattering effects, either directly within a single acoustic propagation computation or via an additional computation. Grid-based methods directly capture scattering effects as part of the acoustic propagation by solving the differential form of a set of propagation equations, such as the linearized Euler equations~\cite{bogey2002computation} and acoustic perturbation equations~\cite{ewert2003acoustic}, with appropriate boundary conditions on scattering surfaces. Grid-based methods accurately capture the effects of complex surface geometries, source representations, and background flows. However, they require a computationally expensive volume discretization, which may be prohibitive in many applications. They also require specialized numerical schemes~\cite{tam1993dispersion} to avoid numerical dispersion and dissipation artifacts~\cite{tam1992discretization} and do not automatically satisfy the Sommerfeld radiation condition~\cite{hu2001stable}. Geometrical acoustics (GA) methods~\cite{thomas2022systematic} approximate acoustic waves as particles that propagate along rays, with models for diffraction along smooth surfaces and at sharp edges. This approximation is formally valid in the high frequency limit but loses accuracy at lower frequencies. While GA methods are highly efficient for high frequency scattering problems, they require the evaluation of ray paths connecting source and observer regions, which may be challenging for moving or distributed sources. Boundary element methods (BEM) and equivalent source methods (ESM) solve integral formulations of a wave equation in the time~\cite{Hu2017TDBEM, Lee2010ESM} or frequency domain~\cite{BARBARINO2018585, tinetti2011fast}. They are typically applied as a separate computation, taking an incident acoustic field computed by an integral acoustic analogy as input to determine a scattered acoustic field that enforces acoustic boundary conditions on scattering surfaces. The scattered field is then added to the incident field to obtain the total acoustic field that includes scattering effects. BEMs solve for unknown layer potential distributions on the scattering surfaces, while ESMs represent the scattered field using equivalent sources located inside a scattering body. This simplifies the numerical implementation by avoiding evaluation of integrals with singular integrands, at the cost of degraded performance for thin bodies and near edges, additional free parameters (including source positions) that must be tuned to different geometries, and reliance on ad hoc stabilization methods~\cite{lee2017use}.

Propellers and rotors are particularly significant sources of aerodynamic noise in a variety of emerging propulsion technologies, including distributed electric propulsion systems~\cite{greenwood2022challenges} and open rotor engines~\cite{brouwer2016scattering}. Propeller blades may be particularly challenging for many scattering prediction approaches, as they constitute rotating, distributed, broadband noise sources with significant low frequency content. Especially for early design analysis, time-domain boundary element methods (TDBEM) have several advantages for sources with these characteristics. Unlike frequency-domain methods, time-domain methods can simulate all frequencies in a broadband signal in a single run~\cite{hu2016assessment}. They can also simulate transient or aperiodic problems and naturally couple with time-domain source data from CFD solvers. BEMs provide a combination of robustness for complex geometries and accuracy and computational efficiency at moderate frequencies.

Despite the advantages of TDBEMs, implementation challenges have limited their widespread adoption. One of the main difficulties is ensuring stability in time integration~\cite{costabel2004time}. Previous TDBEM approaches for aeroacoustic problems have used point collocation discretization strategies~\cite{Hu2017TDBEM, lee2014derivation}, which suffer from instabilities due to non-uniqueness of the discretized system near critical frequencies that correspond to interior resonances of the scattering geometry. Collocation approaches are typically stabilized by applying the Burton-Miller reformulation procedure~\cite{burton1971application}. However, this technique introduces a numerical coupling parameter that must be chosen to achieve stability while obtaining accurate results~\cite{chappell2009choice, zheng2015burton, marburg2016burton, kreuzer2024burton}. Alternatively, full space-time Galerkin discretization has been shown to result in an unconditionally stable formulation~\cite{bamberger1986formulation} that offers a number of numerical and practical advantages. First, the Galerkin approach avoids the introduction of numerical parameters that must be tuned based on the geometry and source characteristics, as is the case for Burton-Miller reformulation. This is critical for use as a predictive tool, where the properties of the solution are not known a priori for each case. Unlike the Burton-Miller collocation approach which can only simulate closed scattering bodies, the space-time Galerkin approach also allows simulation of thin scattering surfaces that do not enclose volumes, which can result in an improved representation and reduced problem size for plate-like scattering geometries. Additionally, Galerkin methods benefit from a quasi-best approximation property in the energy norm that tends to result in more accurate solutions with a given resolution compared to collocation methods~\cite{yu2010development}, and allow relatively straightforward residual-based a posteriori error estimation~\cite{kita2001error,gimperlein2020residual}. The practical use of the space-time Galerkin approach in acoustic scattering problems is currently limited by the substantial cost and complexity of performing a double integration in space and time.

In this paper, we develop a space-time Galerkin time-domain boundary element method for aeroacoustic scattering problems, based on the mathematical formulation of Ha-Duong~\cite{ha2003retarded} for the wave equation. For aeroacoustic applications, we introduce the effect of a mean background flow via variable transformations. We propose an efficient integration approach that overcomes the numerical difficulties of the Galerkin double integrals through careful decomposition of the integrands and integration domains. We validate the method by considering three test cases that reflect different features of realistic aircraft scattering problems. We compare the numerical results with analytical solutions, and investigate some of the numerical properties of the method. Finally, we apply our TDBEM to evaluate scattering and shielding in a trailing edge-mounted propeller configuration, demonstrating successful coupling with CFD source data for numerical simulation of a practical case.


\section{Theory}
\subsection{Governing equations}
The propagation of acoustic waves in a steady irrotational background flow is described by a linearized potential wave equation~\cite{campos200736}:
\begin{equation}
\frac{D_0}{Dt}\left( \frac{1}{c_0^2} \frac{D_0 \phi}{Dt} \right) - \frac{1}{\rho_0} \nabla\cdot(\rho_0 \nabla \phi) = S(x,t), \label{eqn:lpwe}
\end{equation}
where $\phi$ is the acoustic potential, $S$ is a distribution of acoustic sources, $D_0/Dt = \partial/\partial t + u \cdot \nabla$ is the linearized material derivative, $u$ is the background flow velocity, $\rho_0$ is the background density, and $c_0$ is the speed of sound. The acoustic potential is related to the acoustic pressure and acoustic velocity by:
\begin{subequations}
\begin{align}
p & = -\rho_0 \frac{D_0 \phi}{D t}, \\
v & = \nabla \phi.
\end{align}
\end{subequations}

Let $\phi_i$ denote the incident acoustic field, which solves Eq. (\ref{eqn:lpwe}) and is already known (e.g., from an initial acoustic computation or as an analytical expression). The goal of the scattering solver is to compute the unknown scattered field $\phi_s$ such that the total acoustic field $\phi = \phi_i + \phi_s$ solves Eq. (\ref{eqn:lpwe}) together with a sound-hard boundary condition on a scattering surface, representing the impermeability of the scatterer. The scattered field is therefore the homogeneous ($S = 0$) solution of Eq. (\ref{eqn:lpwe}) with a boundary condition residual exactly canceling that of $\phi_i$ for all time.

With suitable simplifying assumptions on the background flow, a variable transformation $(x, t, \phi, S) \rightarrow (\tilde{x}, \tilde{t}, \tilde{\phi}, \tilde{S})$ can exactly or approximately transform Eq. (\ref{eqn:lpwe}) to the standard wave equation (nondimensionalized such that $\tilde{c}_0 = 1$):
\begin{equation}
    \frac{\partial^2 \tilde{\phi}}{\partial \tilde{t}^2} - \nabla_{\tilde{x}}^2 \tilde{\phi} = \tilde{S}(\tilde{x}, \tilde{t}), \label{eqn:wave}
\end{equation}
where $\tilde{(\cdot)}$ denotes the transformed variables. Such variable transformations include the Prandtl-Glauert-Lorentz (PGL) transformation~\cite{hu2019use} for uniform background flows and Taylor's transformation for low Mach number nonuniform background flows~\cite{taylor1978transformation}, as well as various extensions and compositions~\cite{mancini2016integral, tinetti2005aeroacoustic}. The impermeability condition on the scattering surface requires zero surface-normal velocity from both the background and acoustic flow components. This can be enforced exactly for nonuniform background flows around general scatterers, and the acoustic boundary condition therefore requires zero acoustic potential gradient in the surface-normal direction. For a uniform background flow and a scattering body with finite thickness, the background flow must violate the impermeability condition. As shown by Hu et al. (2019)~\cite{hu2019use}, the correct acoustic boundary condition requires zero energy flux (ZEF) through the surface, which reduces to the zero normal velocity boundary condition where the surface is parallel to the background flow direction. Conveniently, the ZEF boundary condition is equivalent to requiring zero normal acoustic potential gradient after the PGL transform is applied. Therefore, both classes of variable transformation result in the same acoustic boundary condition in the transformed variables:
\begin{equation}
    \hat{n}\cdot \tilde{v} = \hat{n}\cdot \nabla_{\tilde{x}} \tilde{\phi} = \frac{\partial \tilde{\phi}}{\partial n} = 0, \label{eqn:neu_BC}
\end{equation}
where $\hat{n}$ denotes the unit normal vector on the (transformed) scattering surface. We can therefore simplify the formulation of the original scattering problem to find $\tilde{\phi}_s$ such that $\tilde{\phi} = \tilde{\phi}_i + \tilde{\phi}_s$ solves Eqs. (\ref{eqn:wave}) and (\ref{eqn:neu_BC}). We also have significant freedom to choose the most appropriate background flow approximation for each scattering problem without adding substantial implementation effort. Because the flow effects are entirely contained in the transformation $\tilde{(\cdot)}$, all problems that can be approximated with these classes of background flows can be solved with a single algorithmic approach for the transformed problem.

\subsection{Boundary integral formulation}
We now derive a boundary integral equation to obtain the scattered field in the transformed variables. For notational clarity, we drop $\tilde{(\cdot)}$ on the transformed variables. Let $\Gamma$ denote a bounded, orientable scattering surface, which does not necessarily enclose a volume. As $\phi_i$ is assumed to be known, the scattered field $\phi_s$ solves the following Neumann problem:
\begin{subequations}
\begin{align}
\frac{\partial^2 \phi_s}{\partial t^2} - \nabla^2 \phi_s &= 0, \, t \in [0,T], \label{eqn:wave_scatter} \\
\frac{\partial \phi_s}{\partial n} = -\frac{\partial \phi_i}{\partial n} &= g, \, x \in \Gamma, \label{eqn:neumann_data} \\
\phi_s &= 0, \, t=0, \label{eqn:wave_ic1} \\
\frac{\partial \phi_s}{\partial t} &= 0, \, t=0, \label{eqn:wave_ic2} 
\end{align}\label{eqn:neumann_prob}
\end{subequations}
where Eqs. (\ref{eqn:wave_ic1}) and (\ref{eqn:wave_ic2}) are initial conditions imposed for uniqueness. The free-space Green's function for Eq. (\ref{eqn:wave_scatter}) is:
\begin{equation}
    G(x,y,t,\tau) = \frac{\delta\left(t - \tau - |x-y|\right)}{4\pi|x-y|},
\end{equation}
where $\delta$ is the Dirac delta function. We define the single and double layer potential distributions on $\Gamma$, $\sigma = \left[\partial \phi_s /\partial n \right]$ and $\psi = [\phi_s]$, where $[\cdot]$ denotes the jump over $\Gamma$ in the direction of the normal vector. With the aid of the initial conditions, we can define $\phi_s$ in terms of the unknown layer potentials $\sigma$ and $\psi$ using the following representation formula:
\begin{equation}
    \phi_s(x,t) = \int_0^t \int_\Gamma \frac{\partial G(x,y,t,\tau)}{\partial n_y}\psi(y, \tau) - G(x,y,t,\tau) \sigma(y, \tau) \, dy d\tau. \label{eqn:repres}
\end{equation}
If we allow $\Gamma$ to represent a thin geometry that does not enclose a volume, the acoustic boundary condition must be allowed to apply on both the upper and lower sides, so the normal acoustic potential gradient jump over the surface is zero. Therefore, we set $\sigma = 0$. By differentiation of Eq. (\ref{eqn:repres}), we obtain a first kind boundary integral equation for the Neumann problem:
\begin{equation}
    \int_0^t\int_\Gamma \frac{\partial^2 G}{\partial n_x \partial n_y}\psi(y, \tau) \, dyd\tau = \frac{\partial}{\partial n_x}\phi_s(x,t) = g(x,t), \label{eqn:rpbie}
\end{equation}
for $x \in \Gamma$. After solving for the double layer potential distribution $\psi$, the scattered field $\phi_s$ is given by the representation formula Eq. (\ref{eqn:repres}) with $\sigma = 0$. Bamberger and Ha-Duong~\cite{bamberger1986formulation} obtain the following weak formulation of Eq. (\ref{eqn:rpbie}):
\begin{align}
    a(\psi, \eta) & = f(\eta), \label{eqn:weak_form} \\
    a(\psi, \eta) & = \int_0^T \int_0^t \int_{\Gamma \times \Gamma} \frac{\partial^2 G}{\partial n_x \partial n_y}\psi(y, \tau) \dot{\eta}(x,t) \, dy dx d\tau dt, \nonumber \\
    f(\eta) & = \int_0^T \int_{\Gamma} g(x,t) \dot{\eta}(x,t) \, dx dt, \nonumber
\end{align}
for test and trial functions belonging to an appropriate function space~\cite{bamberger1986formulation} and where $\dot{(\cdot)}$ denotes time differentiation. The weak formulation is coercive and the resulting Galerkin discretizations are unconditionally stable~\cite{bamberger1986formulation, ha2003retarded}. The bilinear form contains a hypersingular integrand but can be evaluated as~\cite{ha2003retarded,aimi2013neumann}:
\begin{equation}
    a(\psi, \eta) = \frac{1}{4\pi} \int_0^T \int_{\Gamma \times \Gamma} \left(\frac{n_x \cdot n_y}{|x-y|}\dot{\psi}(y,t-|x-y|)\ddot{\eta}(x,t) - \frac{\operatorname{curl}_\Gamma \psi(y,t-|x-y|) \cdot \operatorname{curl}_\Gamma \dot{\eta}(x,t)}{|x-y|} \right) dydxdt,
\end{equation} \label{eqn:bilinear}
where $\operatorname{curl}_\Gamma$ is the tangential curl operator.

\subsection{Discretization}
We approximate the scattering surface $\Gamma$ with a triangulation $\Gamma^h$ consisting of non-overlapping triangular elements $\Gamma_i$. We introduce a uniform temporal discretization of the interval $[0,T]$ with timesteps $t_k = k \Delta t$, where $\Delta t = T/N_t$ is the timestep size. To represent the double layer potential, we use the classical $P^1$ space of piecewise linear basis functions in both space and time, which is the lowest order class with sufficient regularity. The discretized double layer potential $\psi^h$ can then be written in terms of its values at the nodes of $\Gamma^h$ at each timestep:
\begin{equation}
    \psi^h(x,t) = \sum_{i=1}^{N_s} \sum_{k=1}^{N_t} \psi_{ik} \alpha_i(x) \beta_k(t), \, x \in \Gamma^h, \, t \in [0,T],
\end{equation}
where $\alpha_i \in P^1(\Gamma^h)$ are the spatial basis functions, $\beta_k \in P^1(\{t_n\})$ are the temporal basis functions, $N_s$ is the number of nodes, $N_t$ is the number of timesteps, and the coefficient $\psi_{ik}$ is the double layer potential at node $i$ and time $t_k$. We note that on the edges of a thin surface that does not enclose a volume, a discontinuity in the acoustic field over the scattering surface cannot be physically maintained. Therefore, the double layer potential is set to zero at nodes that are on the boundary of the triangulation. This helps ensure numerical stability and reduce the problem size.

Following Ha-Duong~\cite{ha2003retarded}, we choose the test functions $\dot{\eta}^h = \alpha_i(x) \chi_k(t)$, where $\chi_k$ are constant functions supported on $(t_{k-1},t_k]$. Substituting $\psi^h$ and $\dot{\eta}^h$ into Eq. (\ref{eqn:weak_form}), we obtain a linear system:
\begin{subequations}
\begin{align}
    K_{ijkl} \cdot \psi_{ik} & = f_{jl}, \label{eqn:lin_sys}\\
    K_{ijkl} & = \frac{1}{4\pi} \int_0^T \int_{\Gamma^h \times \Gamma^h}
    \begin{aligned}[t] & \left( \frac{n_x \cdot n_y}{|x-y|}\alpha_i(y) \alpha_j(x) \dot{\beta}_k(t-|x-y|) \dot{\chi}_l(t) \right. \\
    & \,\, {}- \left. \frac{\operatorname{curl}_\Gamma \alpha_i(y) \cdot \operatorname{curl}_\Gamma \alpha_j(x)}{|x-y|} \beta_k(t-|x-y|) \chi_l(t) \right) \, dydxdt, \end{aligned} \label{eqn:infl_coeff_disc} \\
    f_{jl} & = \int_0^T \int_{\Gamma^h} g(x,t) \alpha_j(x) \chi_l(t) \, dxdt. \label{eqn:neumann_data_disc}
\end{align}
\end{subequations}
The influence matrix coefficients $K_{ijkl}$ are nonzero only if $t - |x-y| \in (t_{k-1}, t_{k+1})$ for $t \in (t_{l-1}, t_l]$, $x \in \operatorname{supp}(\alpha_j)$, and $y \in \operatorname{supp}(\alpha_i)$. We define:
\begin{subequations}
    \begin{align}
        d_{ij} & = \min\{|x-y|: x\in \operatorname{supp}(\alpha_j), y \in \operatorname{supp}(\alpha_i)\}, \\
        D_{ij} & = \max\{|x-y|: x\in \operatorname{supp}(\alpha_j), y \in \operatorname{supp}(\alpha_i)\}.
    \end{align}
\end{subequations}
Then the coefficients are nonzero for $t_l - t_k \in (d_{ij} - \Delta t, D_{ij} + 2 \Delta t)$. As the scattering body is assumed to be finite, $D_{ij}$ is bounded and therefore there exists $N_b \in \mathbb{N}$ such that $K_{ijkl} = 0$ for $l-k > N_b$. Importantly, the coefficients also do not depend directly on $k$ and $l$, but only the difference between them. Therefore, Eq. (\ref{eqn:lin_sys}) can be rewritten as:
\begin{equation}
    \sum_{m=0}^{N_b} K_{ij}^{m} \psi_{i,l-m} = f_{jl}, \, K_{ij}^m = K_{ij,1,1+m}.
\end{equation}
This leads to a marching on in time (MOT) algorithm:
\begin{equation}
    K_{ij}^0 \psi_{il} = f_{jl} - \sum_{m=1}^{N_b} K_{ij}^{m} \psi_{i,l-m}, \, l=1,\dots,N_t. \label{eqn:MOT}
\end{equation}
For long simulations, the number of matrix coefficients that must be computed is small and independent of the simulation time, and the MOT scheme only requires the repeated solution of a relatively small $N_s \times N_s$ linear system for each timestep.

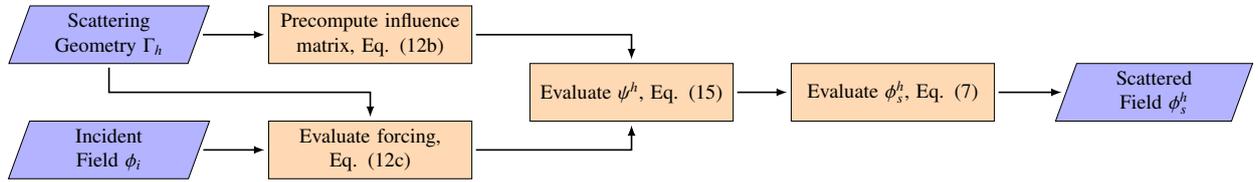
\begin{figure}[h]
\centering
\resizebox{\textwidth}{!}{
\begin{tikzpicture}[node distance=2cm,every node/.style={font=\linespread{1}\selectfont}]
\node (geo) [io] {Scattering\\ Geometry $\Gamma_h$};
\node (infl) [process, right of=geo, xshift=2.5cm] {Precompute influence\\ matrix, Eq. (\ref{eqn:infl_coeff_disc})};

\node (incid) [io, below of=geo] {Incident\\ Field $\phi_i$};
\node (forcing) [process, right of=incid, xshift=2.5cm] {Evaluate forcing,\\ Eq. (\ref{eqn:neumann_data_disc})};

\node (MOT) [process, right of=forcing, xshift=2.5cm, yshift=1cm] {Evaluate $\psi^h$, Eq. (\ref{eqn:MOT})};
\node (prop) [process, right of=MOT, xshift=2.5cm] {Evaluate $\phi^h_s$, Eq. (\ref{eqn:repres})};
\node (out) [io, right of=prop, xshift=2.5cm] {Scattered\\ Field $\phi^h_s$};

\draw[thick,->] (geo) -- (infl);
\draw[thick,->] (incid) -- (forcing);
\draw[thick,->] (geo) -- ++(0,-1) -| (forcing);
\draw[thick,->] (infl) -| (MOT);
\draw[thick,->] (forcing) -| (MOT);
\draw[thick,->] (MOT) -- (prop);
\draw[thick,->] (prop) -- (out);
\end{tikzpicture}}
\caption{Flow chart showing the process of computing the scattered acoustic field from the scattering geometry and incident field. The incident field is provided through coupling with another solver (e.g., CFD with integral-based acoustic analogy) or as an analytical expression.}
\label{fig:flow_chart}
\end{figure}

The complete numerical solution procedure for the Neumann problem Eq. (\ref{eqn:neumann_prob}) is summarized in Fig. \ref{fig:flow_chart}. The influence matrix is evaluated directly from the scattering geometry, and can therefore be precomputed. The forcing $f_{jl}$ is evaluated from the incident field and the scattering geometry. The double layer potential $\psi^h$ is computed with the MOT algorithm using the influence matrix and forcing. Finally, the scattered field is obtained from the double layer potential via the representation formula. These three steps can either be performed sequentially or combined into a single time-marching algorithm to accumulate the scattered field.

The required spatial and temporal resolutions are not limited by stability requirements, but to avoid information loss, the node spacing and timestep must be smaller than half the wavelength and period, respectively, of the highest resolved frequency due to the Nyquist sampling limit. The node spacing and timestep can be related by a CFL parameter:
\begin{equation}
    \mbox{CFL} = \frac{\Delta t}{\Delta x},
\end{equation}
where typically $\mbox{CFL} \leq 1$ to satisfy causality~\cite{ha2003retarded}, although $\mbox{CFL} > 1$ is still stable. For a surface mesh with approximately uniform element sizes, the node spacing can be estimated as $\Delta x \approx \sqrt{2A}$, where $A$ is the average element area.

\subsection{Numerical challenges}
The main numerical difficulty is in the evaluation of the double space-time integral to compute the influence matrix coefficients. In particular, the integrand contains several different singularities that degrade the accuracy of standard quadrature rules. We define:
\begin{subequations}
    \begin{align}
        f^m_1(R) & = \int_{0}^T \dot{\beta}_{1}(t-R) \dot{\chi}_{1+m}(t) \, dt, \\
        f^m_2(R) & = \int_{0}^T {\beta}_{1}(t-R) {\chi}_{1+m}(t) \, dt,
    \end{align}
\end{subequations}
where $\beta_k(t)$ and $\chi_k(t)$ are the temporal basis functions for the trial and test functions. The functions $f^m_1$ and $f^m_2$ are piecewise polynomials and have discontinuities or discontinuous derivatives at $R=n\Delta t$ for $n\in[m-2,m+1]$. The influence matrix coefficients can be rewritten by interchanging the order of integration:
\begin{equation}
    K_{ij}^m = \frac{1}{4\pi} \int_{\Gamma^h \times \Gamma^h} \left( \frac{n_x \cdot n_y}{|x-y|}\alpha_i(y) \alpha_j(x) f_1^m(|x-y|) - \frac{\operatorname{curl}_\Gamma \alpha_i(y) \cdot \operatorname{curl}_\Gamma \alpha_j(x)}{|x-y|} f_2^m(|x-y|) \right) \, dydx.
\end{equation}
Two types of singularities are apparent: an $\mathcal{O}(|x-y|^{-1})$ singularity due to the Green's function, and light cone singularities at $|x-y|=n\Delta t$ due to the discontinuities in $f^m_1$ and $f^m_2$. Building on previous work by Aimi et al.~\cite{aimi2013neumann}, we consider an iterated approach, taking the inner and outer spatial integrals sequentially to enable geometric decomposition of the integration domain around the singularities.

\begin{figure}[h]
\centering
\begin{subfigure}[b]{0.3\textwidth}
\centering
\begin{tikzpicture}[scale=1.5]
\begin{scope}
\def\angEl{35} 
\def\angAz{-100} 

\tikzset{xyplane/.style={cm={cos(\angAz),sin(\angAz)*sin(\angEl),-sin(\angAz),
                          cos(\angAz)*sin(\angEl),(0,0)}}}

\tikzset{vert1/.style={cm={0,1,-sin(\angAz),cos(\angAz)*sin(\angEl),(0,0)}}}

\fill[gray!30,xyplane] (0.35,0.45) -- (1.3,0.7) -- (0.3,1.5);
\draw[xyplane] (0.35,0.45) -- (1.3,0.7) -- (0.3,1.5) -- cycle;
\draw[xyplane,->] (1,1.2) -- (0.65,0.85) node[below right, pos=0] {$\Gamma_s$};

\draw[vert1] (0.6,0.6) -- (1.1,-0.2) -- (0.3,-0.5) -- cycle;
\draw[vert1,->] (1,0.5) -- (0.75,0.1) node[right, pos=0, inner sep=2pt] {$\Gamma_t$};

\draw[vert1] (0,0) -- (0.44,0);
\draw[vert1,->] (0.7,0) -- (1.5,0) node[above, inner sep=2pt] {$z'$};

\draw[xyplane,->] (0,0) -- (2,0) node[pos=1, below] {$x'$};
\draw[xyplane,->] (0,0) -- (0,1.8) node[pos=1, right] {$y'$};

\draw[xyplane,dashed,->] (0,0) -- (1.88,0.68) node[pos=0.9, right, inner sep=2pt] {$r'$};
\draw[xyplane,dashed,<->] (1.6,0) arc (0:20:1.6) node[pos=0.5, below] {$\theta'$};

\fill[black] (0,0.7) circle (1pt) node[left, inner sep=1pt] {$x_t$};
\end{scope}
\end{tikzpicture}
\caption{Coordinates for inner integration}
\label{fig:inner_integal_coords}
\end{subfigure}
\begin{subfigure}[b]{0.3\textwidth}
\centering
\begin{tikzpicture}[scale=2.25]
\begin{scope}
\draw[line join=round] (0.3,1.5) -- (0.35,0.45) -- (1.3,0.7) -- cycle;

\draw[dashed,pattern={Lines[angle=-45, distance=4pt]}, pattern color=orange,line join=round] (0,0) -- (1.3,0.7) -- (0.3,1.5) -- cycle;
\draw[dashed,pattern={Lines[angle=45, distance=4pt]}, pattern color=blue,line join=round] (0,0) -- (0.35,0.45) -- (1.3,0.7) -- cycle;
\draw[dashed,pattern={Lines[angle=90, distance=4pt]}, pattern color=teal,line join=round] (0,0) -- (0.35,0.45) -- (0.3,1.5) -- cycle;

\draw[->] (1,0.35) -- (0.5,0.4) node[right, pos=0, color=blue, inner sep=2pt] {$\Gamma_s^1$};
\draw[->] (1.05,1.25) -- (0.45,0.65) node[above right, pos=0, color=orange!80!black] {$\Gamma_s^2$};
\draw[->] (0.15,1.5) -- (0.22,0.6) node[above, pos=0, color=teal!80!black, inner sep=2pt] {$\Gamma_s^3$};

\draw[->] (0,0) -- (1.6,0) node[right] {$x'$}; 
\draw[->] (0,0) -- (0,1.85) node[above] {$y'$};
\end{scope}
\end{tikzpicture}
\caption{Decomposition of $\Gamma_s$ into component triangles}
\label{fig:inner_integal_decomp}
\end{subfigure}
\begin{subfigure}[b]{0.3\textwidth}
\centering
\begin{tikzpicture}[scale=2.25]
\begin{scope}
\draw (0,0) -- (1.3,-0.6) -- (1.3,0.7) -- cycle;
\draw[dashed] (1.35,0) arc (0:50:1.35);
\draw[dashed,<->] (0,0) -- (0.95,0.95) node[above left, pos=0.55] {$r'_n$};
\draw[dashed] (1.35,0) arc (0:-30:1.35);
\draw[dashed] (1.3,-0.9) -- (1.3,-0.6);
\draw[dashed] (1.3,1.035) -- (1.3,0.7);
\draw[dashed,<->] (0,-0.8) -- (1.3, -0.8) node[above, pos=0.5] {$x''_0$};

\draw[pattern={Lines[angle=45, distance=4pt]}, pattern color=gray,line join=round] (0,0) -- (1.226,-0.566) -- (1.226,-0.566) arc (-24.8:-15.6:1.35) -- cycle;
\draw[pattern={Lines[angle=90, distance=4pt]}, pattern color=gray,line join=round] (0,0) -- (1.3,-0.364) -- (1.3,0.364) -- cycle;
\draw[pattern={Lines[angle=-45, distance=4pt]}, pattern color=gray,line join=round] (0,0) -- (1.189,0.64) -- (1.189,0.64) arc (28.3:15.6:1.35) -- cycle;

\draw[->] (0,0) -- (1.6,0) node[right] {$x''$}; 
\draw[<->] (0,-0.93) -- (0,0.93) node[above] {$y''$};
\end{scope}
\end{tikzpicture}
\caption{Component triangle decomposition}
\label{fig:inner_integal_light_cone}
\end{subfigure}
\caption{Geometry for inner spatial integration: (a) definition of local Cartesian ($x'$, $y'$, $z'$) and cylindrical ($r'$, $\theta'$, $z'$) coordinate systems; (b) decomposition of $\Gamma_s$ into three component triangles with a vertex at the origin; (c) decomposition of the domain formed by the intersection of a component triangle and circle into sector and triangle regions.}
\label{fig:inner_integal_geo}
\end{figure}
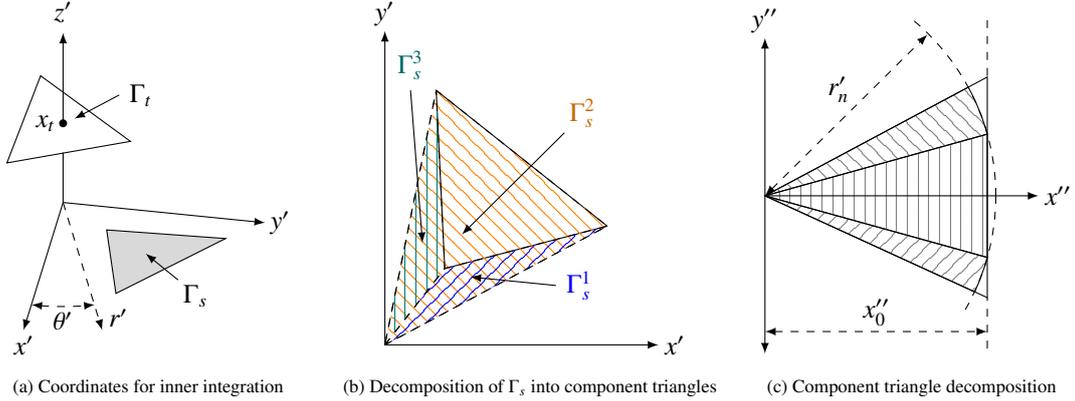

The matrix coefficient $K^m_{ij}$ is assembled from integrals over pairs of triangles associated with the nodes $i$ and $j$. We therefore focus on the integral over one pair of triangles $\Gamma_s \in \operatorname{supp}(\alpha_i)$ and $\Gamma_t \in \operatorname{supp}(\alpha_j)$. To define the inner integral, let us fix a target point $x_t \in \Gamma_t$. We introduce a local Cartesian coordinate system $(x',y',z')$ such that the \mbox{$x'y'$-plane} coincides with the plane of the triangle $\Gamma_s$ and the origin lies at the projection of the target point, as depicted in Fig. \ref{fig:inner_integal_coords}. We define $r' = \sqrt{x'^2 + y'^2}$ and $\theta' = \arctan(y'/x')$. These variables are convenient for the inner integration because $|x-y|=\sqrt{r'^2 + z'^2}$, so the Jacobian of the transformation to the cylindrical coordinates $(r', \theta', z')$ analytically cancels the $\mathcal{O}(|x-y|^{-1})$ singularity at the origin, while the functions $f^m_1$ and $f^m_2$ are independent of $\theta'$. Therefore, the light cone singularities occur at constant integration limits in $r'$. The integrand expressed in cylindrical coordinates is now regular on the subdomains obtained by taking intersections of $\Gamma_s$ with annuli defined by successive light cone boundaries $r'_n = \sqrt{n^2\Delta t^2 - z'^2}$, although these subdomains are complicated in the general case.

To simplify the subdomains, we note that the annulus between successive light cone boundaries is the difference between two circular domains, and that $\Gamma_s$ can be decomposed as the union or difference of three component triangles, each with a vertex at the origin (see Fig. \ref{fig:inner_integal_decomp}), resulting in constant integration limits in $\theta'$. The inner integral can then be expressed as the sum of several regular integrals on the domains formed by the intersections between the component triangles and circles with radii $r'_n$. Let $(x'',y'',z')$ and $(r',\theta'',z')$ denote an in-plane rotation of the coordinate systems $(x',y',z')$ and $(r',\theta',z')$ such that the edge of the component triangle opposite the origin is parallel to the $y''$-axis. As illustrated in Fig. \ref{fig:inner_integal_light_cone}, the intersection of a component triangle with a circle of radius $r'_n$ is in general the union of two sectors and triangle, each of which are simple and easy to define.

Existing integration strategies use 2D quadrature rules to take the regular integrals, which have the form:
\begin{equation}
    I_P = \int_{\theta''_{min}}^{\theta''_{max}} \int_0^{r'(\theta)} P(R(r)) \alpha_j(r,\theta) \frac{r}{R(r)} \, dr d\theta, \label{eqn:aimi_integral}
\end{equation}
where $r'(\theta'') = r'_n$ for a sector and $x_0''/\cos(\theta'')$ for a triangle, $R(r) = \sqrt{r^2+z'^2}$, and $P(R)$ is a polynomial. The spatial basis function $\alpha_j$ can be expressed in these coordinates as:
\begin{equation}
    \alpha_j(r',\theta'') = c_1 + r'(c_2 \cos(\theta'') + c_3 \sin(\theta'')).
\end{equation}
To reduce the computational cost of applying a 2D quadrature rule to the inner integrals, which must subsequently be integrated again over $\Gamma_t$, we introduce an additional decomposition and analytical integration of Eq. (\ref{eqn:aimi_integral}) that results in a 1D numerical integral. We note that $P(R)$ can be decomposed into monomial terms containing integer powers of $R$. Therefore, we are interested in the following 2D integrals:
\begin{subequations}
\begin{align}
    I^c_{\alpha} & = \int_{\theta''_{min}}^{\theta''_{max}} I_{\alpha,1}(r'(\theta),z') \, d\theta, \\
    I^{x''}_{\alpha} & = \int_{\theta''_{min}}^{\theta''_{max}} I_{\alpha,2}(r'(\theta),z') \cos{\theta} \, d\theta, \\
    I^{y''}_{\alpha} & = \int_{\theta''_{min}}^{\theta''_{max}} I_{\alpha,2}(r'(\theta),z') \sin{\theta} \,  d\theta,
\end{align} \label{eqn:2D_int}
\end{subequations}
where:
\begin{equation}
    I_{\alpha,\beta}(r,z) = \int_0^r \left(\rho^2 + z^2\right)^{\frac{\alpha}{2}} \rho^\beta \, d\rho.
\end{equation}
These are standard integrals with simple exact expressions (see \ref{app:integrals}). For our choice of temporal basis functions, we implement the exact expressions for the parameter ranges $\alpha = -1,0,1$ and $\beta=1,2$. For the sector subdomains, the radial and angular integrations are separable and the entire 2D integrals in Eq. (\ref{eqn:2D_int}) are easily taken exactly. For the triangle subdomains, exact expressions for Eq. (\ref{eqn:2D_int}) can also be obtained~\cite{carley2013analytical, milroy1997elastostatic}. However, except when very high accuracy is required, the exact integrals are less numerically stable and more expensive to evaluate then applying numerical quadrature~\cite{milroy1997elastostatic}. We therefore use a 1D Gauss-Legendre quadrature rule to perform the angular integration on the triangle subdomain. Compared to existing approaches that rely on 2D quadrature schemes, the computational cost is substantially reduced by the exact radial integration and the corresponding quadrature dimension reduction.

The outer spatial integral over $\Gamma_t$ is subsequently taken with standard Gaussian quadrature schemes, such as those obtained by Dunavant~\cite{dunavant1985high}. Although not completely smooth~\cite{ostermann2010numerical}, the outer integrand is sufficiently regular to provide good convergence without resorting to a specialized quadrature strategy.

The decomposition-based quadrature strategy developed for the inner spatial integration is also used to evaluate the scattered field from the layer potential with the representation formula. Substituting $\psi^h$ into Eq. (\ref{eqn:repres}), we obtain:
\begin{equation}
    \phi^h_s(x,t) = \sum_{i=1}^{N_s} \sum_{k=1}^{N_t} \frac{\psi_{ik}}{4\pi} \int_{\Gamma^h} \frac{\hat{n}_y \cdot (x-y)}{|x-y|}\left(\frac{\beta_k(t-|x-y|)}{|x-y|^2} + \frac{\dot{\beta}_k(t-|x-y|)}{|x-y|}\right)  \alpha_i(y) \, dy.
\end{equation}
The temporal basis function $\beta_k(t)$ is piecewise polynomial, as is its derivative, resulting in light cone singularities at $|x-y| = t - n\Delta t$ for $n \in [k-1, k+1]$. In the local cylindrical coordinates $(r', \theta', z')$ defined for each $\Gamma_i \in \Gamma^h$, $\hat{n}_y \cdot (x-y) = z'$ and is independent of $r'$ and $\theta'$. The integral on each $\Gamma_i$ can therefore be written as a combination of the integrals in Eq. (\ref{eqn:2D_int}) using the same geometric and algebraic decompositions presented above for the Galerkin inner integrals.

\section{Validation results}
To validate the TDBEM formulation and implementation, we consider three test cases with features that are representative of realistic aeroacoustic scattering problems and compare the numerical predictions with analytical reference solutions. The key characteristics of the test cases are summarized in Table \ref{tab:validation_summary}.

\begin{table}[h]
\centering
\caption{Summary of geometry, flow, and source attributes for the validation cases.}
\begin{tabular}{llll}
\hline
Case & Geometric features & Background flow & Source type \\
\hline
Sphere & Smooth curvature, enclosed body & None & Harmonic, broadband \\
Disk & Flat surface, sharp edges & None & Harmonic \\
Plane & Flat surface & Uniform & Transient \\
\hline
\end{tabular}
\label{tab:validation_summary}
\end{table}

In all three cases, we consider a homogeneous medium and define nondimensional coordinates such that $\rho_0 = 1$ and $c_0 = 1$. We prescribe an incident field generated by a monopole point source, given by an analytical expression~\cite{rienstra2004introduction}:
\begin{equation}
    \phi_i(x,t) = -\frac{Q(t - R)}{4\pi R^*}, \label{eqn:monopole_src}
\end{equation}
where $R = (R^* - |x-y|M_r)/\beta^2$, $R^* = |x-y|\sqrt{M_r^2 + \beta^2}$, $M_r = M_\infty \cdot (x-y)/|x-y|$, $\beta^2 = 1-|M_\infty|^2$, $M_\infty$ is the uniform freestream Mach vector, $y$ is the monopole source location, and $Q$ is the monopole source strength function.

\subsection{Sphere}
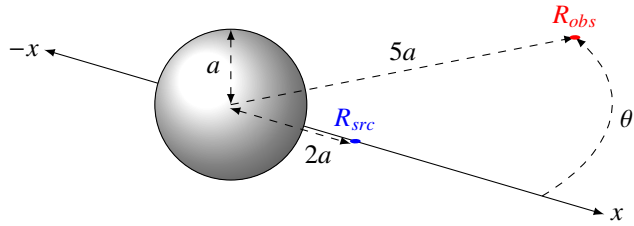
\begin{figure}[h]
\centering
\begin{tikzpicture}

\def\R{1} 
\def\angEl{25} 
\def\angAz{-125} 

\tikzset{xyplane/.style={cm={cos(\angAz),sin(\angAz)*sin(\angEl),-sin(\angAz),
                              cos(\angAz)*sin(\angEl),(0,0)}}}

\filldraw[ball color=white] (0,0) circle (\R);

\draw[xyplane,<-] (0,-3*\R) node[left] {$-x$} -- (0,-1.17*\R);
\draw[xyplane,->] (0,1.17*\R) -- (0,6*\R) node[right] {$x$};

\draw[xyplane,->,dashed] (0,0) -- (-1.73*2.5*\R,2.5*\R) node[pos=0.5,above,inner sep=2pt] {$5a$};
\draw[xyplane,->,dashed] (0,5*\R) arc (90:150:5*\R) node[pos=0.5,right,inner sep=2pt] {$\theta$};

\draw[<->,dashed] (0,0) -- (0,\R) node[pos=0.5,left,inner sep=2pt] {$a$};

\draw[xyplane,<->,dashed] (0.1,0.05*\R) -- (0.1,2*\R) node[pos=0.73,below,inner sep=2pt] {$2a$};
\fill[xyplane,blue] (0,2*\R) circle (2pt) node[above,inner sep=2pt] {$R_{src}$};
\fill[xyplane,red] (-1.73*2.5*\R,2.5*\R) circle (2pt) node[above,inner sep=2pt] {$R_{obs}$};
\end{tikzpicture}
\caption{Sphere geometry, showing the source (blue dot, $R_{src}$) and observer (red dot, $R_{obs}$) positions.}
\label{fig:sphere_geo}
\end{figure}

The first validation case is acoustic scattering by a sphere in a medium at rest. This case is representative of a smoothly curved, finite scattering surface that encloses a volume. Because the background flow is at rest, no variable transformation is applied. The geometry consists of a sphere of radius $a$ centered on the origin, with a point source placed at a radius of $2a$ along the positive $x$-axis, and observers located in a circular array centered on the origin at a radius of $5a$ (see Fig. \ref{fig:sphere_geo}). Unless otherwise noted, the sphere is discretized with $5120$ faces and $2562$ nodes, and computations are performed at $\mbox{CFL} = 0.5$.

\begin{figure}[h]
    \centering
    \includegraphics[trim={0 0 0 615},clip,width=0.6\textwidth]{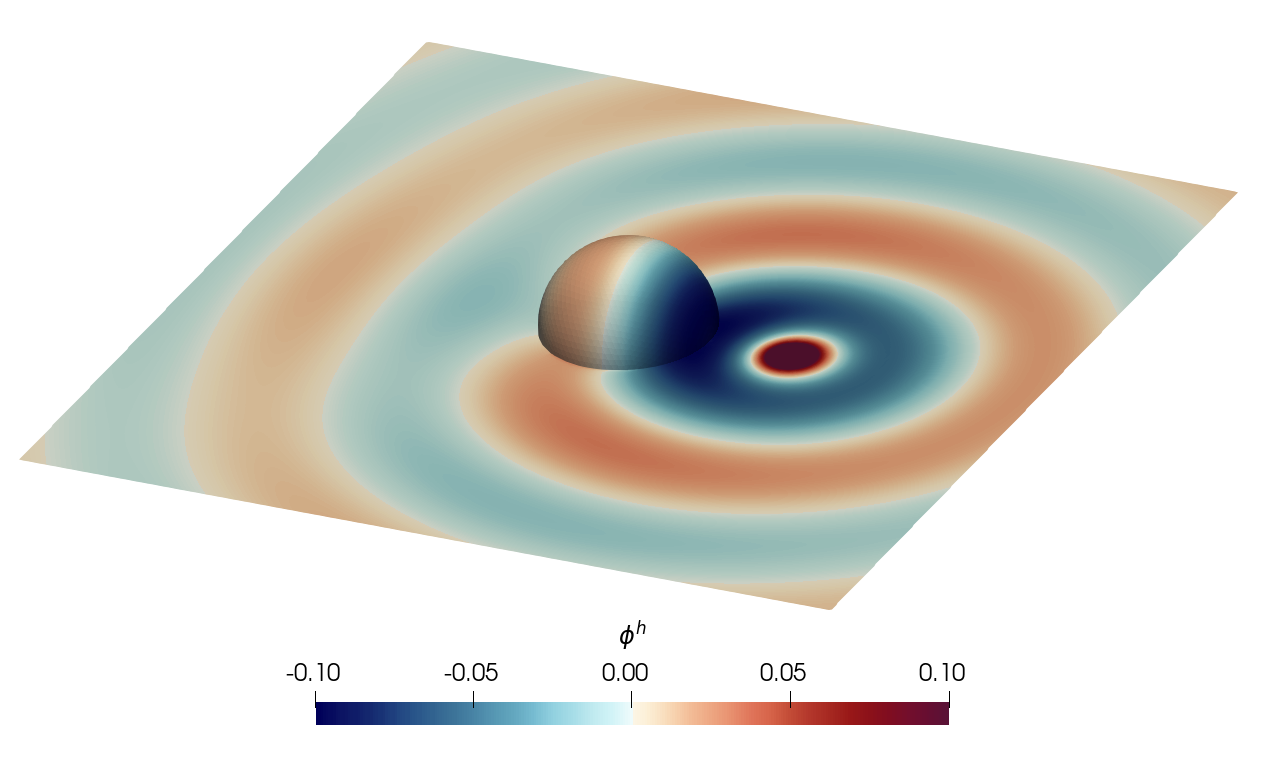}\\
    \begin{subfigure}[b]{0.45\textwidth}
         \centering
         \includegraphics[trim={0 150 0 40},clip,width=\textwidth]{sphere_snapshot_k2_perspective.png}
         \caption{$ka=2$}
         \label{fig:sphere_harmonic_snapshot_k2}
     \end{subfigure}
    \begin{subfigure}[b]{0.45\textwidth}
         \centering
         \includegraphics[trim={0 150 0 40},clip,width=\textwidth]{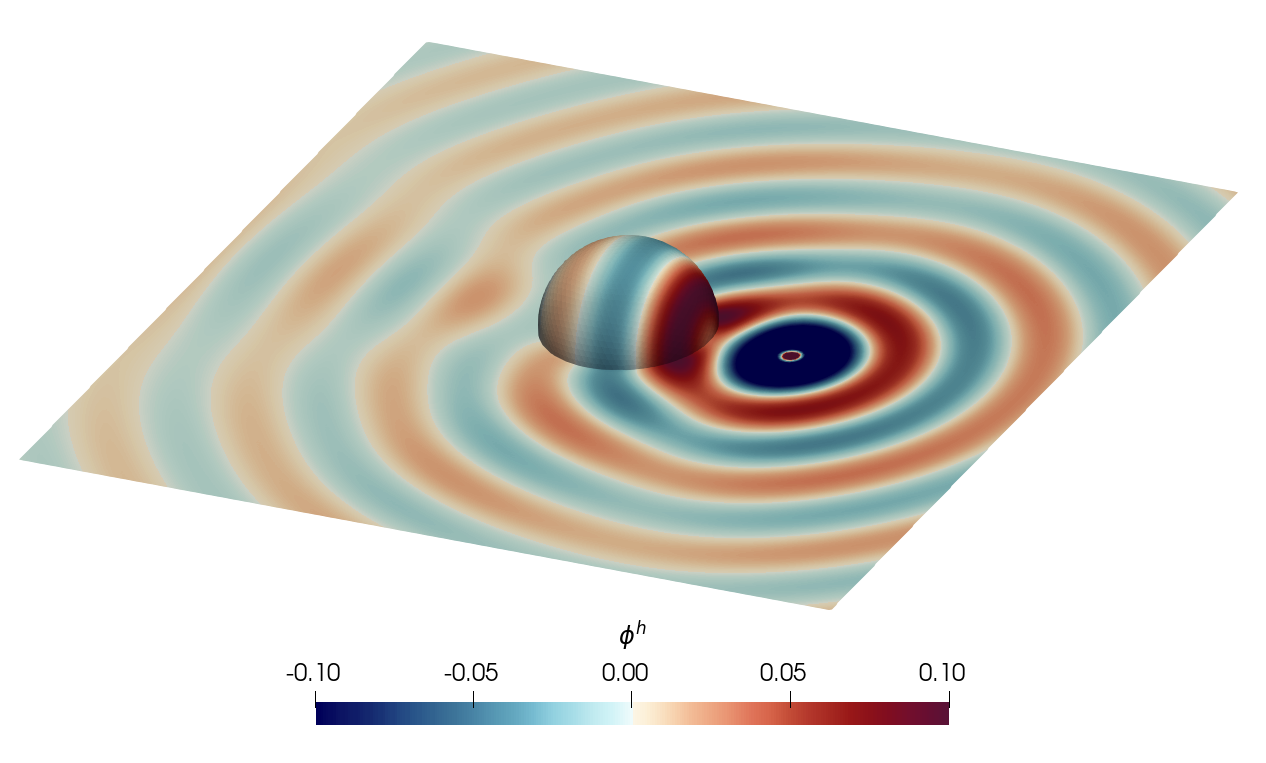}
         \caption{$ka=4$}
         \label{fig:sphere_harmonic_snapshot_k4}
     \end{subfigure}\\
    \begin{subfigure}[b]{0.45\textwidth}
         \centering
         \includegraphics[trim={0 150 0 40},clip,width=\textwidth]{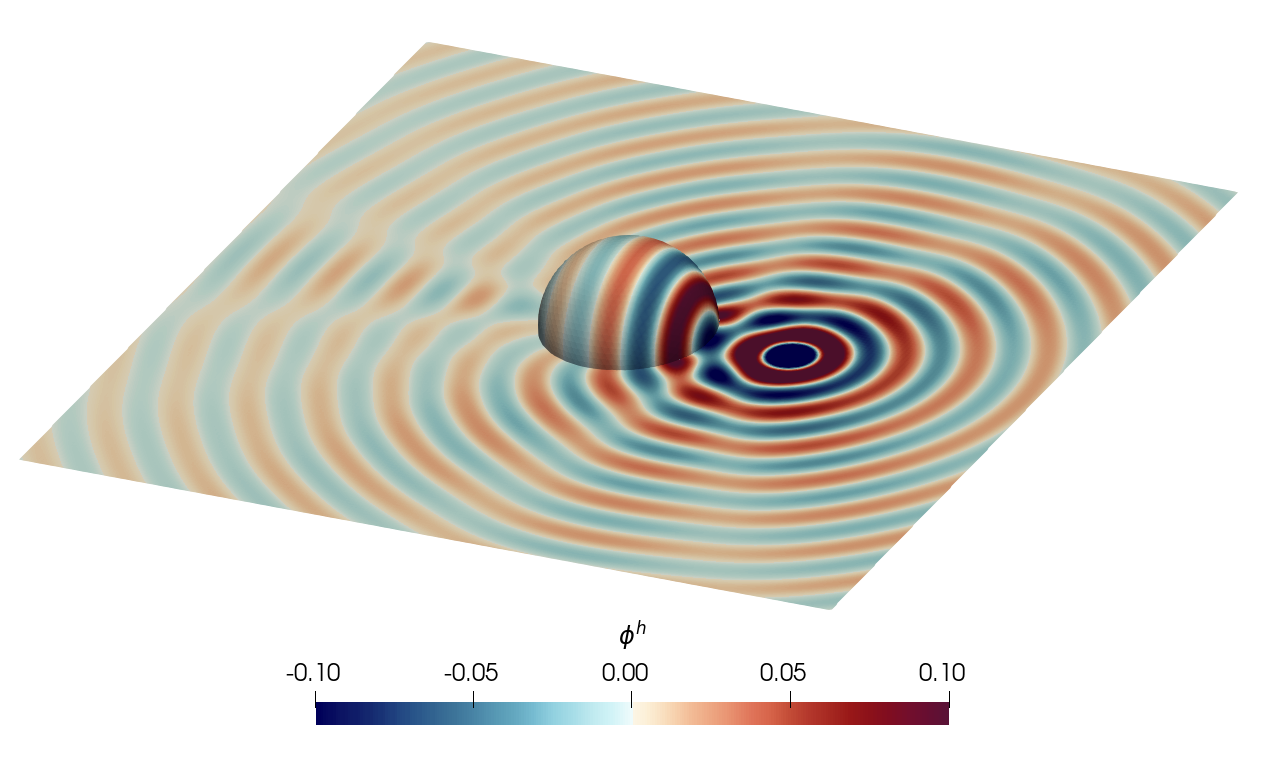}
         \caption{$ka=8$}
         \label{fig:sphere_harmonic_snapshot_k8}
     \end{subfigure}
    \begin{subfigure}[b]{0.45\textwidth}
         \centering
         \includegraphics[trim={0 150 0 40},clip,width=\textwidth]{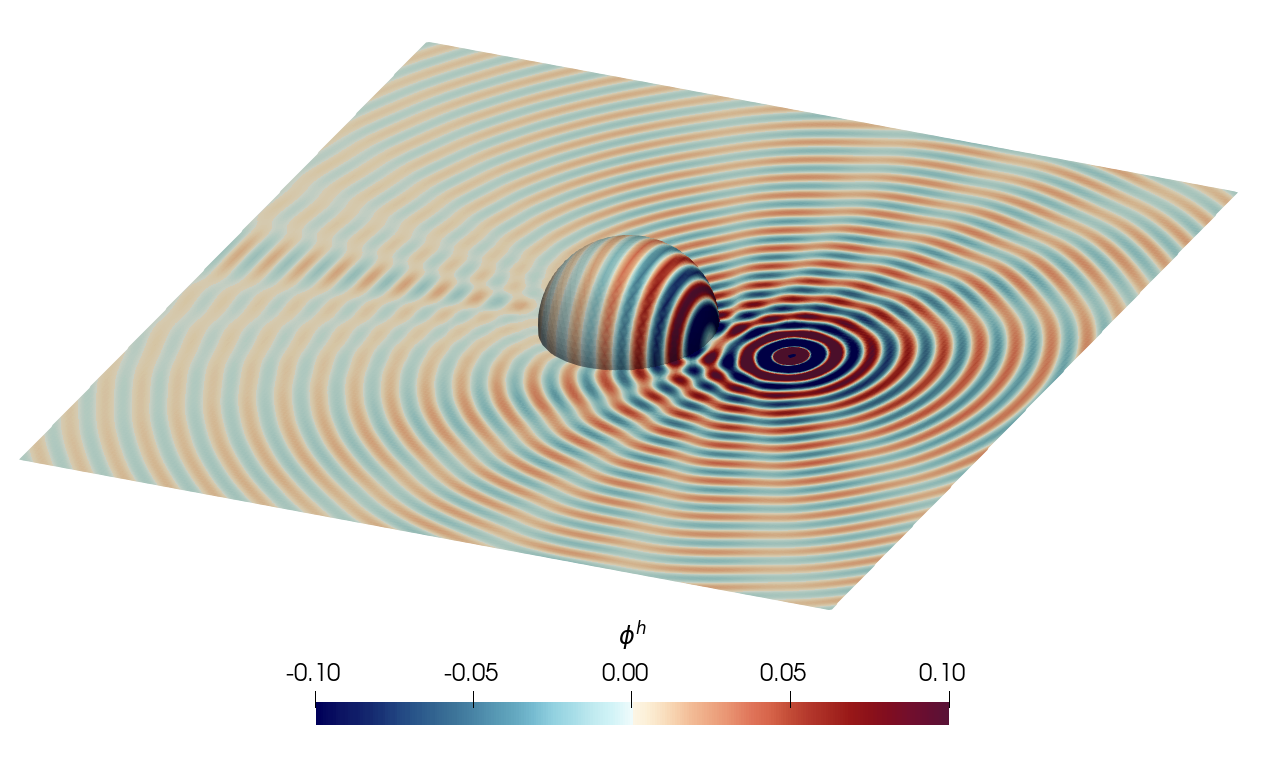}
         \caption{$ka=16$}
         \label{fig:sphere_harmonic_snapshot_k16}
     \end{subfigure}
    \caption{Visualizations of the predicted acoustic field of a harmonic point source scattered by a sphere at selected frequencies: (a) $ka=2$; (b) $ka=4$; (c) $ka=8$; (d) $ka=16$. The instantaneous acoustic potential is shown at nondimensional time $t=50$.}
    \label{fig:sphere_harmonic_snapshots}
\end{figure}

\begin{figure}[h]
    \centering
    \begin{subfigure}[b]{0.45\textwidth}
         \centering
         \includegraphics[trim={0 3 0 3},clip,width=\textwidth]{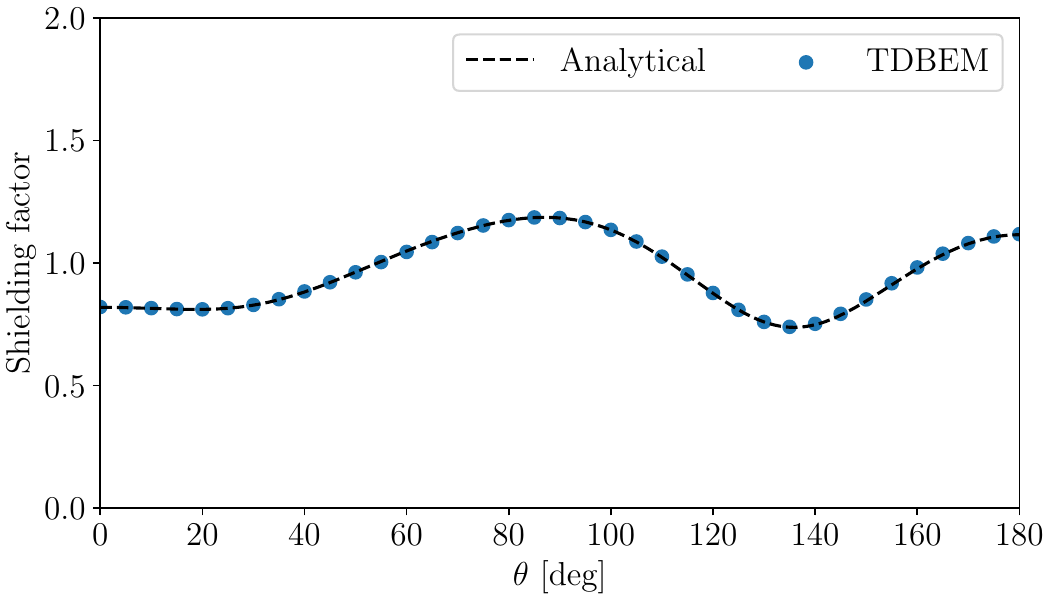}
         \caption{$ka=2$}
         \label{fig:sphere_harmonic_directivity_k2}
     \end{subfigure}
    \begin{subfigure}[b]{0.45\textwidth}
         \centering
         \includegraphics[trim={0 3 0 3},clip,width=\textwidth]{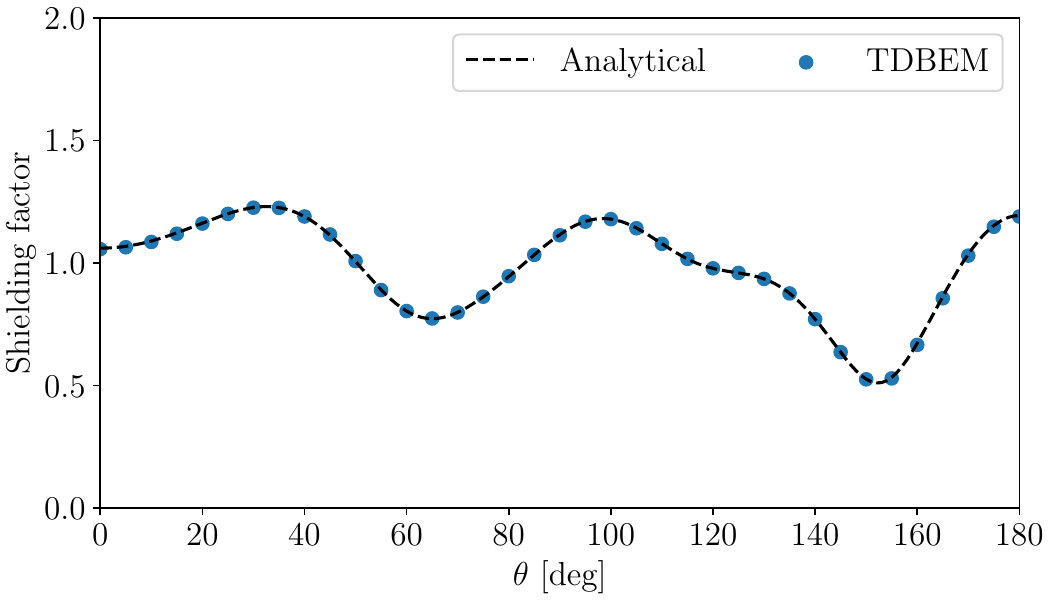}
         \caption{$ka=4$}
         \label{fig:sphere_harmonic_directivity_k4}
     \end{subfigure}\\
    \begin{subfigure}[b]{0.45\textwidth}
         \centering
         \includegraphics[trim={0 3 0 3},clip,width=\textwidth]{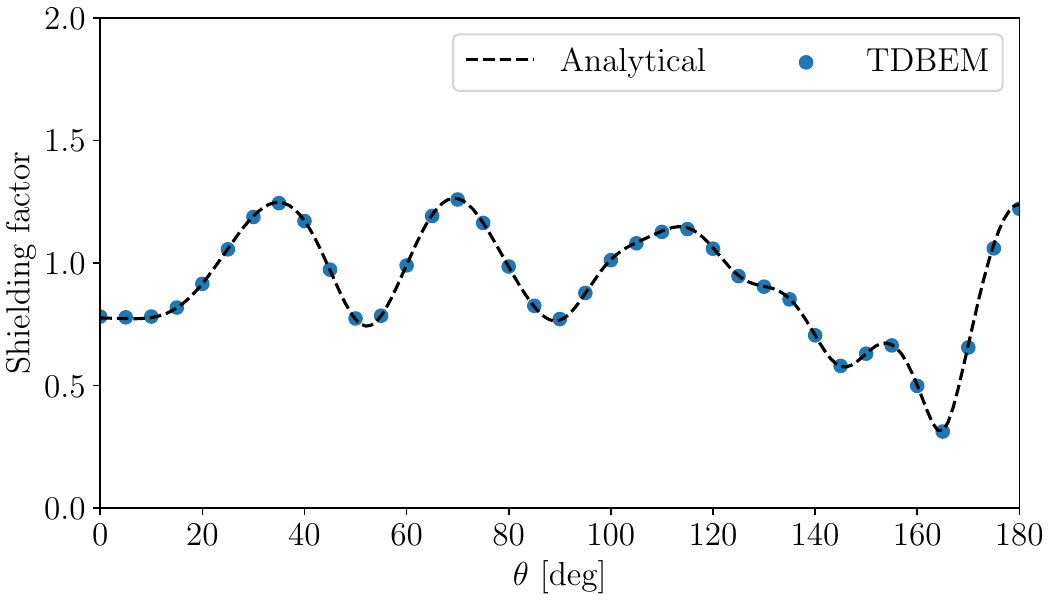}
         \caption{$ka=8$}
         \label{fig:sphere_harmonic_directivity_k8}
     \end{subfigure}
    \begin{subfigure}[b]{0.45\textwidth}
         \centering
         \includegraphics[trim={0 3 0 3},clip,width=\textwidth]{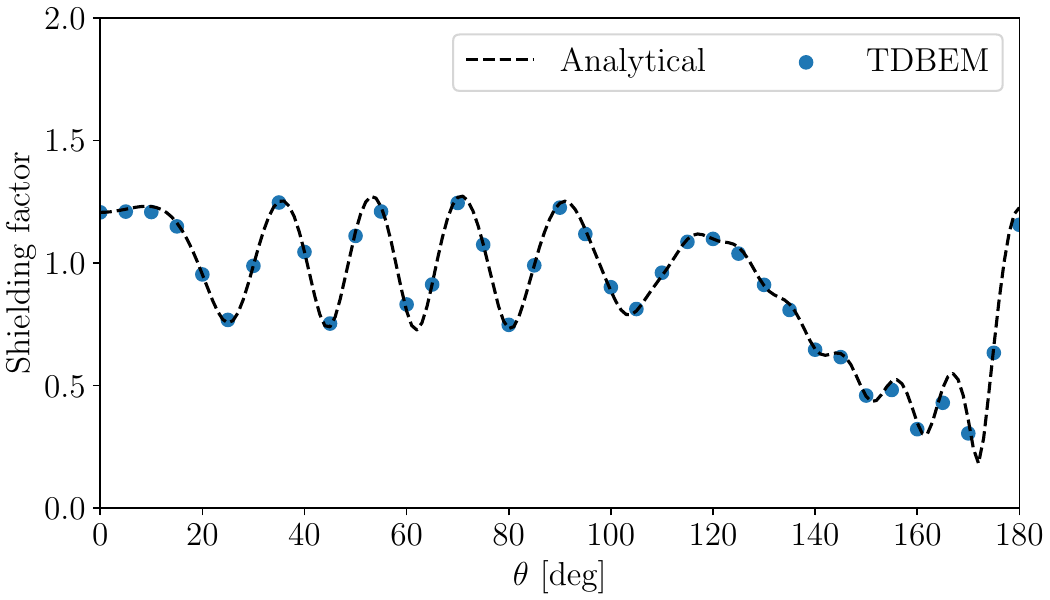}
         \caption{$ka=16$}
         \label{fig:sphere_harmonic_directivity_k16}
     \end{subfigure}
    \caption{Comparison of predicted and analytical shielding factors for scattering of a harmonic point source by a sphere at selected frequencies: (a) $ka=2$; (b) $ka=4$; (c) $ka=8$; (d) $ka=16$.}
    \label{fig:sphere_harmonic}
\end{figure}

We first consider forcing at a constant harmonic frequency, corresponding to $Q(t) = -\cos(kt)$. The analytical solution for this problem is given by Bowman et al.~\cite{bowman1969electromagnetic}. Because the numerical solution is initialized from a quiescent initial condition, the numerical data are collected well after the start of the simulation when initial transient effects have exited the relevant parts of the domain. The acoustic potential distributions $\phi^h$ for harmonic forcing at $ka = 2,4,8,16$ are shown in Fig. \ref{fig:sphere_harmonic_snapshots} on the sphere surface and on an observer plane $x \in [-5,5]$, $y \in [-5,5]$. The visualizations show the increased influence of the scattering body and complexity of the total field with increasing frequency.

To assess the scattering effects, we compute the shielding factor $\gamma_t$, which is the ratio of the amplitudes of the total and incident fields:
\begin{equation}
    \gamma_t = \frac{\phi_{rms}}{\phi_{i,rms}},
\end{equation}
where $(\cdot)_{rms}$ denotes the root mean square. The computed shielding factor directivities on the circular observer array are compared with the analytical solution for $ka = 2,4,8,16$ in Fig. \ref{fig:sphere_harmonic}. At all frequencies and polar angles, the agreement between the predicted and analytical shielding factors is strong. At lower frequencies, the solutions are essentially coincident, while some numerical error becomes visible as the frequency increases because the wavelength and period decrease relative to the fixed spatial and temporal resolution.

\begin{figure}[h]
    \centering
    \begin{subfigure}[b]{0.45\textwidth}
         \centering
         \includegraphics[trim={0 40 0 40},clip,width=\textwidth]{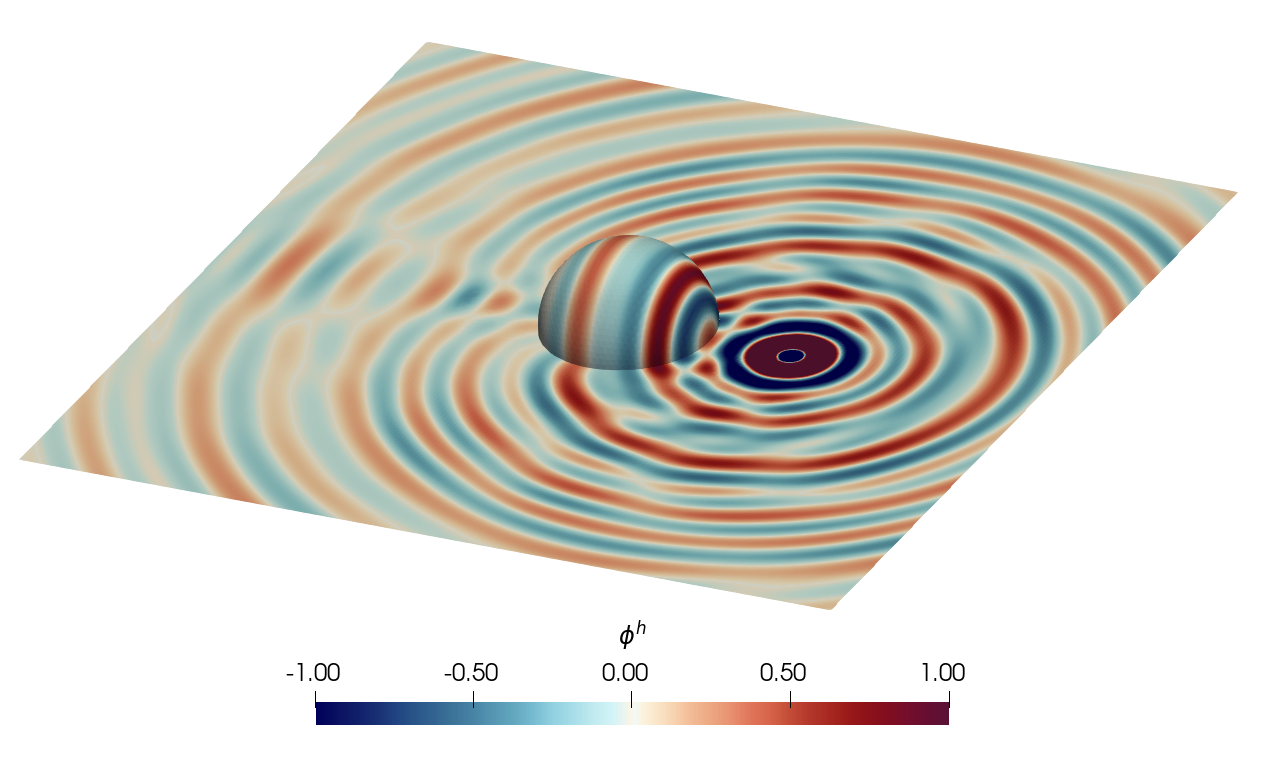}
         \caption{Instantaneous acoustic potential at $t=50$}
         \label{fig:sphere_bb_snapshot}
     \end{subfigure}
    \begin{subfigure}[b]{0.45\textwidth}
         \centering
         \includegraphics[width=\textwidth]{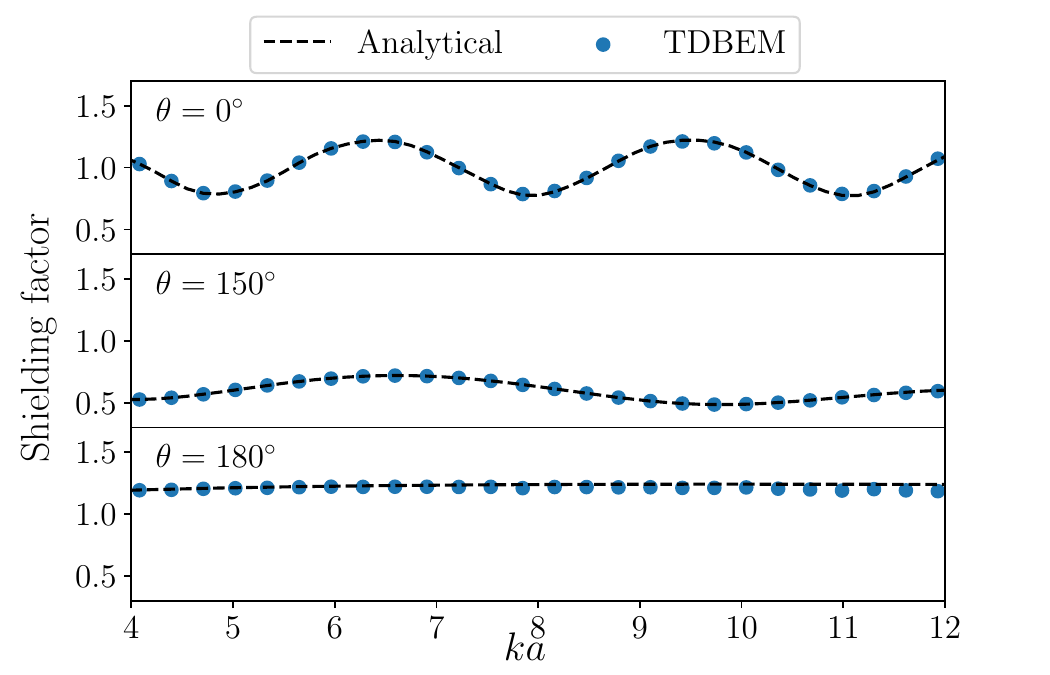}
         \caption{Shielding factor spectra at $\theta=0^\circ, 150^\circ, 180^\circ$}
         \label{fig:sphere_bb_slices}
     \end{subfigure}\\
    \begin{subfigure}[b]{0.45\textwidth}
         \centering
         \includegraphics[width=\textwidth]{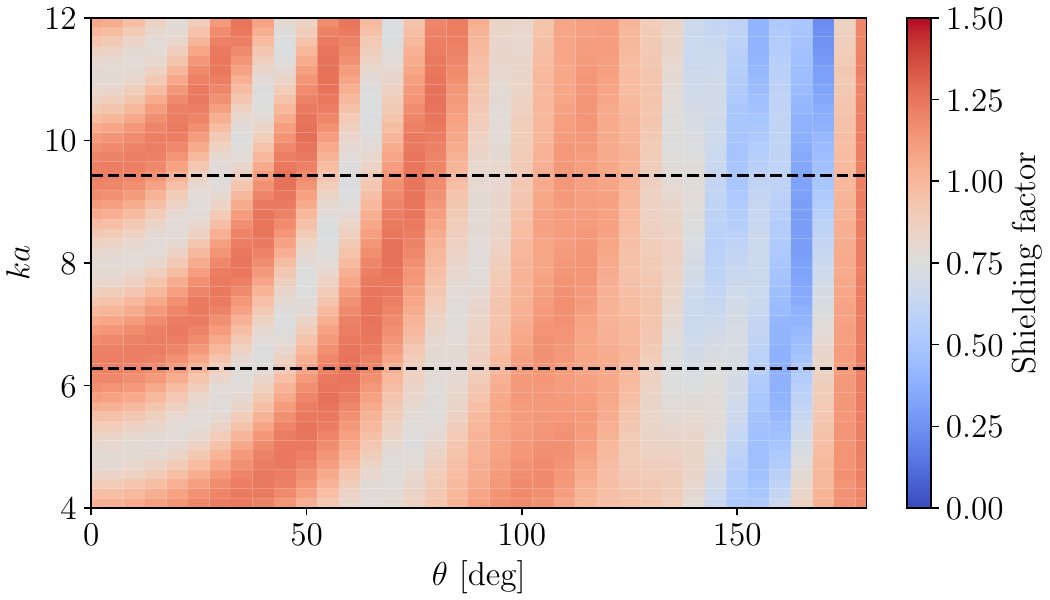}
         \caption{Shielding factor (numerical prediction, $\gamma_t^h$)}
         \label{fig:sphere_bb_total}
     \end{subfigure}
    \begin{subfigure}[b]{0.45\textwidth}
         \centering
         \includegraphics[width=\textwidth]{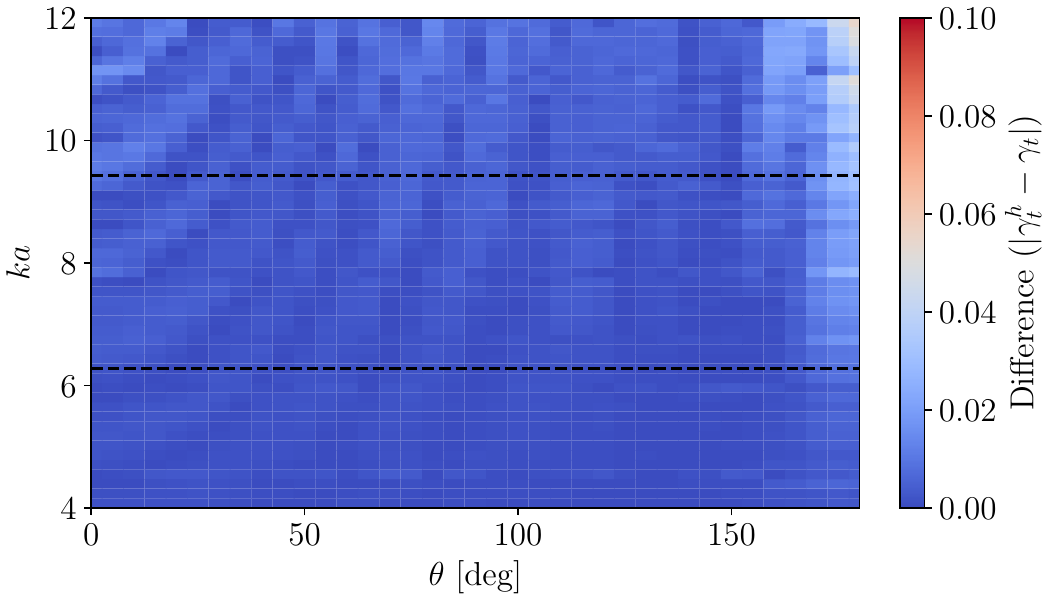}
         \caption{Shielding factor prediction error $|\gamma_t^h - \gamma_t|$}
         \label{fig:sphere_bb_error}
     \end{subfigure}
    \caption{Comparison of predicted and analytical scattering of a broadband point source by a sphere: (a) visualization of predicted acoustic potential on sphere surface and observer plane at nondimensional time $t = 50$; (b) comparison of predicted and analytical shielding factor spectra at selected polar angles; (c) spectrogram of the predicted shielding factor over the circular observer array; (d) spectrogram of the shielding factor prediction error $|\gamma_t^h - \gamma_t|$ over the circular observer array.}
    \label{fig:sphere_bb}
\end{figure}

To establish the performance of the scattering method for broadband signals, we apply band-limited white noise forcing. The signal is synthesized in the frequency domain by assigning uniform amplitudes and random phases to Fourier coefficients within the range $[k_{min},k_{max}]$. The time-domain source strength is then obtained by applying the inverse Fourier transform. To construct the analytical solution, the harmonic solution provided by Bowman et al.~\cite{bowman1969electromagnetic} is evaluated at each frequency and multiplied by the corresponding Fourier coefficient. The predicted acoustic field for a band-limited signal between $k_{min}a = 4$ and $k_{max}a = 12$ is shown in Fig. \ref{fig:sphere_bb_snapshot}. The predicted and analytical shielding factor spectra are then compared at selected polar angles on the circular observer array in Fig. \ref{fig:sphere_bb_slices}, showing very close agreement. The predicted shielding factor and prediction error on the full observer array are presented as spectrograms in Figs. \ref{fig:sphere_bb_total} and \ref{fig:sphere_bb_error}. The prediction error is small at all frequencies and polar angles but is noticeably increased near the Poisson spot at $\theta=180^\circ$ (directly behind the obstacle). The shielding factor error is increased in this region because the magnitude of the scattered field (which contains the prediction error) is large compared to the magnitude of the incident field, not because the prediction of the scattered field itself is less accurate at these angles. The error also increases with frequency due to the fixed resolution. The dashed horizontal lines in Figs. \ref{fig:sphere_bb_total} and \ref{fig:sphere_bb_error} indicate critical frequencies $ka = n\pi$, which are associated with interior resonances that cause instabilities and require stabilization in point collocation TDBEM and ESM approaches. The present method simulates these critical frequencies without instability or increased prediction error.

\begin{figure}[h]
    \centering
    \begin{subfigure}[b]{0.45\textwidth}
         \centering
         \includegraphics[trim={0 3 0 2},clip,width=\textwidth]{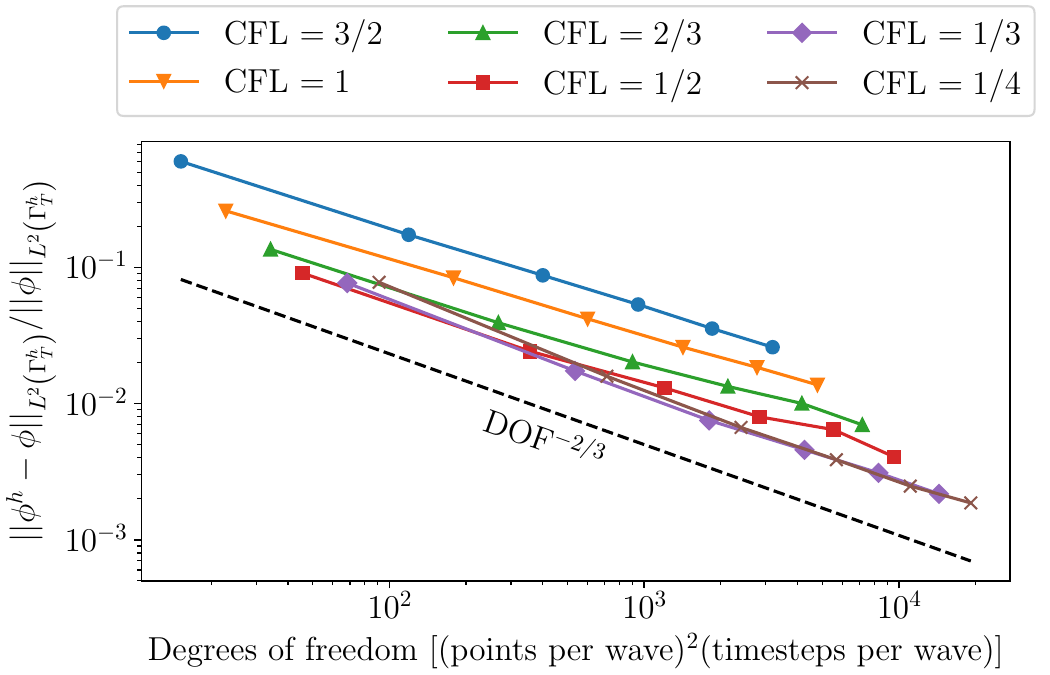}
         \caption{Joint spatial and temporal refinement at fixed $\mbox{CFL}$}
         \label{fig:sphere_conv_line}
     \end{subfigure}
    \begin{subfigure}[b]{0.46\textwidth}
         \centering
         \includegraphics[trim={0 3 0 4},clip,width=\textwidth]{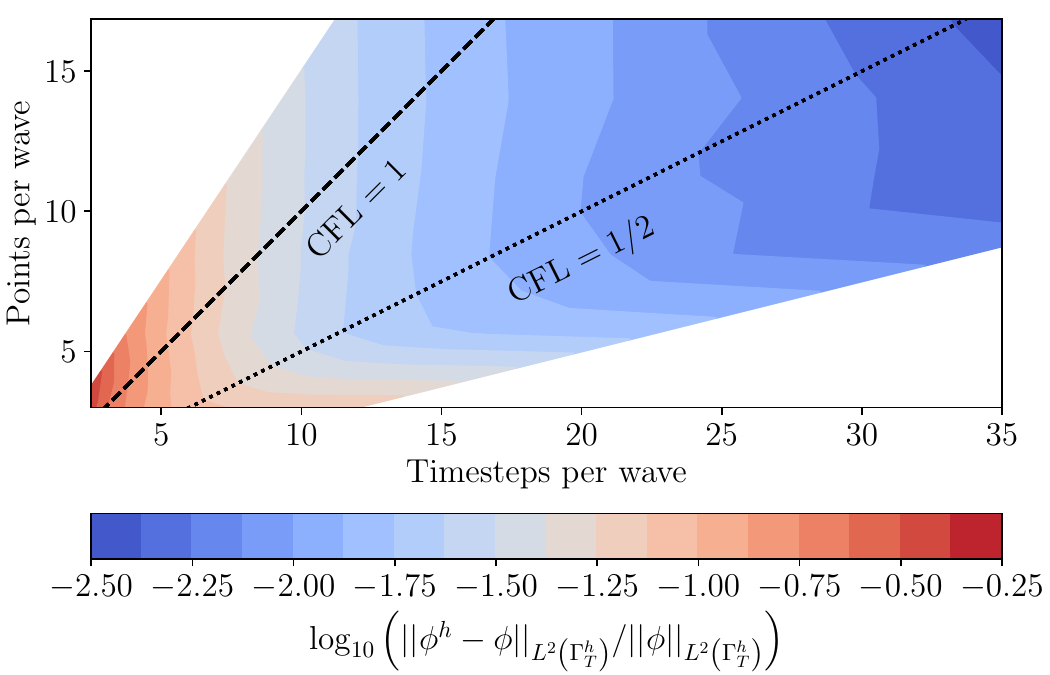}
         \caption{Separate spatial and temporal refinement}
         \label{fig:sphere_conv_cntr}
     \end{subfigure}
    \caption{Refinement study of sphere scattering at $ka = 8$: (a) convergence at fixed $\mbox{CFL}$, showing second order convergence in space and time; (b) convergence with separate spatial and temporal refinement, where fixed $\mbox{CFL} = 1$ and $\mbox{CFL} = 0.5$ are emphasized with dashed and dotted lines, respectively.}
    \label{fig:sphere_conv}
\end{figure}

Finally, to investigate the resolution requirements and convergence properties, we perform a refinement study at the harmonic frequency $ka=8$ across several meshes ranging from $320$ faces and $162$ nodes to $11520$ faces and $5762$ nodes with varying CFL parameters. The error is measured by taking the $L^2$ norm of the instantaneous acoustic potential error $\phi^h - \phi$ on the scattering surface $\Gamma^h$ over the transient-free portion of the simulation time $t \in [t_{min}, T]$. The spatial and temporal resolutions are measured relative to the wavelength $\lambda$ and period $T_p$ by the number of points per wave ($\lambda/\Delta x$) and timesteps per wave ($T_p/\Delta t$), respectively. The second order convergence rate in space and time at fixed $\mbox{CFL}$ is demonstrated in Fig. \ref{fig:sphere_conv_line}, where the relative error is plotted against the relative number of space-time degrees of freedom (DOF) used to describe the layer potential (defined such that $\mbox{degrees of freedom} = (\mbox{points per wave})^2 (\mbox{timesteps per wave})$). The theoretical second order convergence rate in both $\Delta x$ and $\Delta t$ corresponds to the $\mbox{DOF}^{-2/3}$ reference line, which closely matches the observed convergence rates at all $\mbox{CFL}$ values. To explore the significance of the $\mbox{CFL}$ parameter, the error under separate spatial and temporal refinement is presented as a contour plot in Fig. \ref{fig:sphere_conv_cntr}. The error contours have a distinctive ``L'' shape that strongly supports joint spatial and temporal refinement at $\mbox{CFL} \approx 0.5$. For a fixed spatial or temporal resolution, performing temporal or spatial refinement, respectively, results in a rapid initial error reduction towards $\mbox{CFL} \approx 0.5$ before the error stagnates upon further refinement. This suggests that the total error is dominated by temporal error for $\mbox{CFL} \gg 0.5$ and spatial error for $\mbox{CFL} \ll 0.5$. The $\mbox{CFL}$ value does not affect stability, and no change in the convergence behavior is observed for $\mbox{CFL} > 1$.

\FloatBarrier

\subsection{Disk}
\begin{figure}[h]
\centering
\begin{tikzpicture}

\def\R{1} 

\def\angEl{25} 
\def\angAz{-125} 

\tikzset{xyplane/.style={cm={cos(\angAz),sin(\angAz)*sin(\angEl),-sin(\angAz),
                              cos(\angAz)*sin(\angEl),(0,0)}}}
\tikzset{yzplane/.style={cm={cos(\angAz),sin(\angAz)*sin(\angEl),
                              0,cos(\angEl),(0,0)}}}

\draw[yzplane, black, fill=gray!40] circle[radius=\R];

\draw[xyplane,<-] (0,-3*\R) node[left] {$-x$} -- (0,-0.6*\R);
\draw[xyplane,->] (0,0) -- (0,6*\R) node[right] {$x$};

\draw[xyplane,->,dashed] (0,0) -- (-1.73*2.5*\R,2.5*\R) node[pos=0.5,above,inner sep=2pt] {$5a$};
\draw[xyplane,->,dashed] (0,5*\R) arc (90:150:5*\R) node[pos=0.5,right,inner sep=2pt] {$\theta$};

\draw[yzplane,<->,dashed] (0,0) -- (0,\R) node[pos=0.55,right,inner sep=2pt] {$a$};

\draw[xyplane,<->,dashed] (0.1,0.05*\R) -- (0.1,2*\R) node[pos=0.5,below,inner sep=2pt] {$2a$};
\fill[xyplane,blue] (0,2*\R) circle (2pt) node[above,inner sep=2pt] {$R_{src}$};
\fill[xyplane,red] (-1.73*2.5*\R,2.5*\R) circle (2pt) node[above,inner sep=2pt] {$R_{obs}$};

\end{tikzpicture}
\caption{Disk geometry, showing the source (blue dot, $R_{src}$) and observer (red dot, $R_{obs}$) locations.}
\label{fig:disk_geo}
\end{figure}
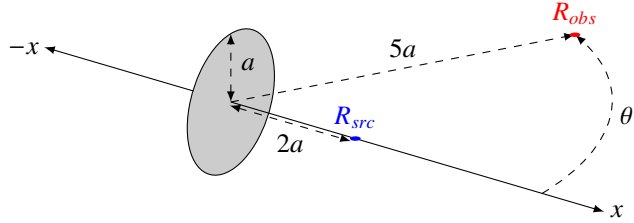

The second validation case is acoustic scattering by a circular disk in a medium at rest. This case is representative of a flat, finite scattering surface with sharp edges. We represent the disk as a thin surface that does not enclose a volume, using a single layer of elements. The ability to represent thin surfaces with a single layer is important because it greatly reduces the computational problem size for applicable cases. Because the background flow is at rest, no variable transformation is applied. The scatterer is a disk of radius $a$ centered on the origin and oriented normal to the $x$-axis. A point source is placed at a radius of $2a$ along the positive $x$-axis, and observers are located in a circular array centered on the origin at a radius of $5a$ (see Fig. \ref{fig:disk_geo}). The disk is discretized with $2400$ faces and $1261$ nodes, and computations are performed at $\mbox{CFL} = 0.5$.

\begin{figure}[h]
    \centering
    \includegraphics[trim={0 0 0 615},clip,width=0.6\textwidth]{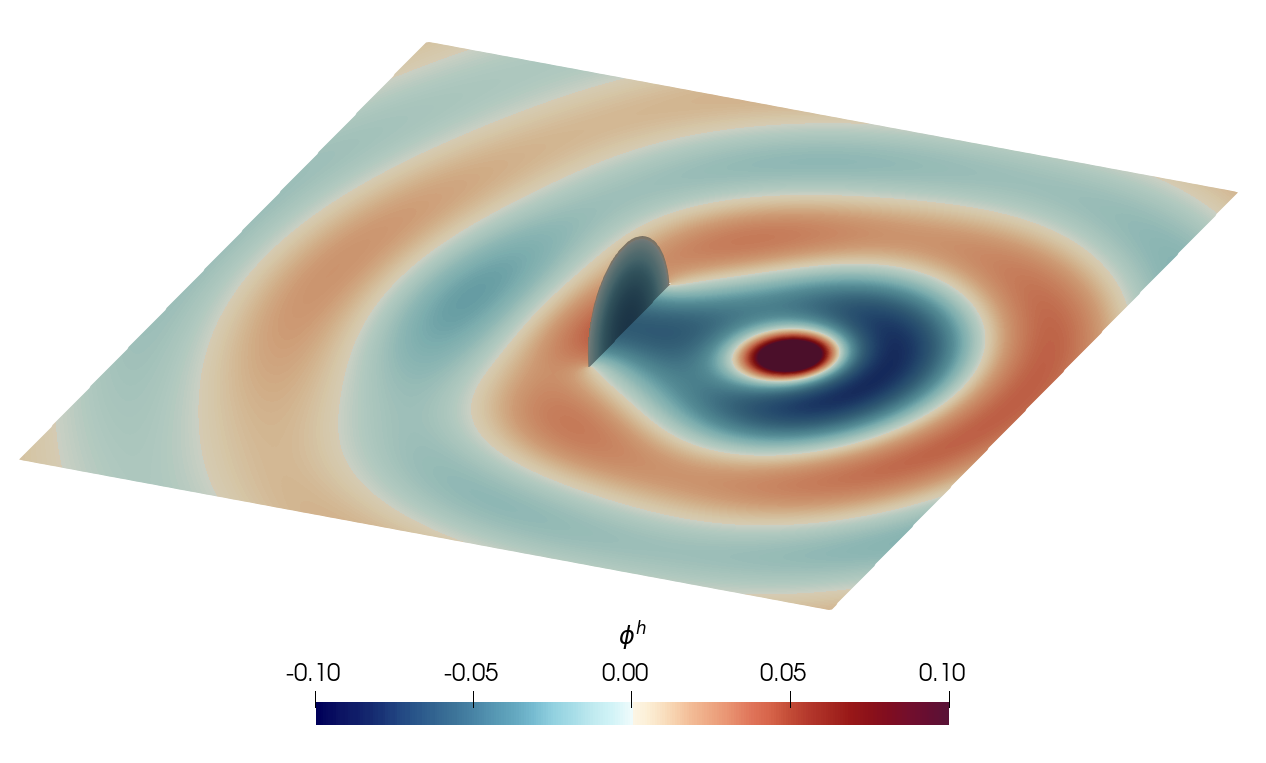}\\
    \begin{subfigure}[b]{0.45\textwidth}
         \centering
         \includegraphics[trim={0 150 0 40},clip,width=\textwidth]{disk_k2_perspective.png}
         \caption{$ka=2$.}
         \label{fig:disk_harmonic_snapshot_k2}
     \end{subfigure}
    \begin{subfigure}[b]{0.45\textwidth}
         \centering
         \includegraphics[trim={0 150 0 40},clip,width=\textwidth]{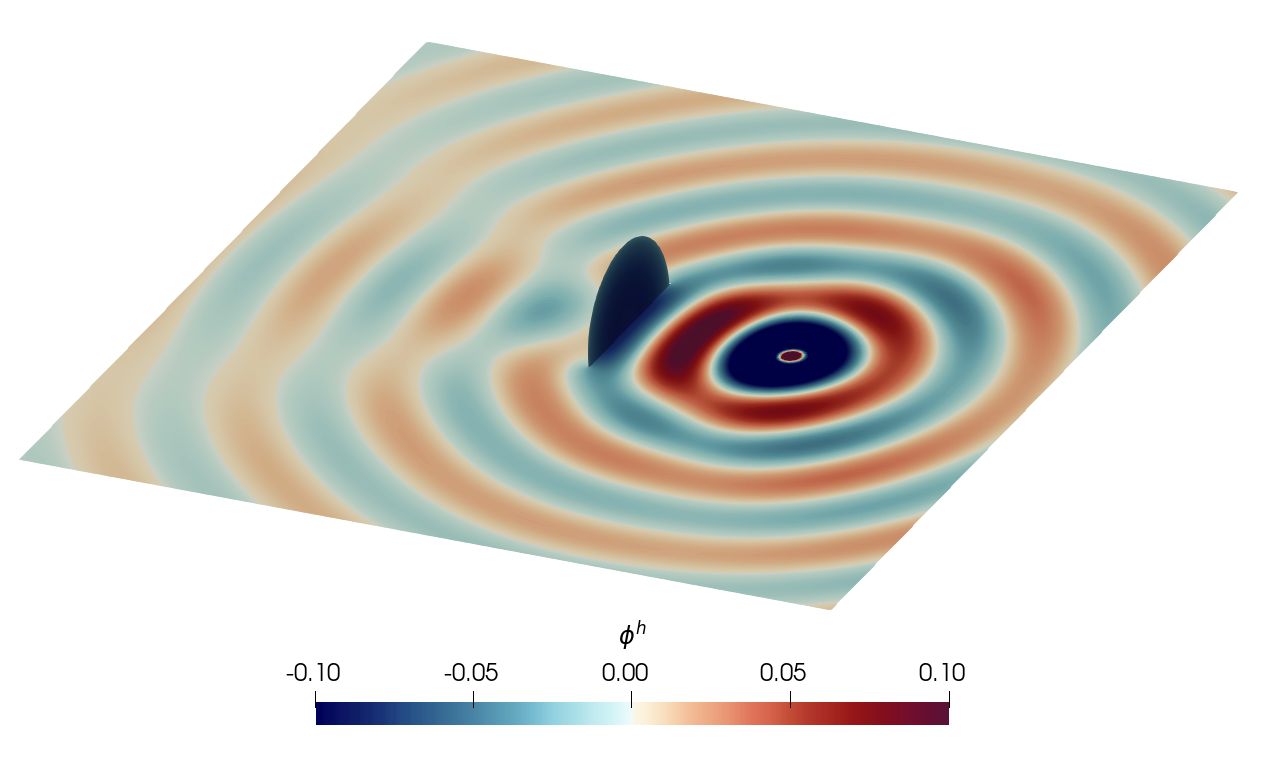}
         \caption{$ka=4$.}
         \label{fig:disk_harmonic_snapshot_k4}
     \end{subfigure}\\
    \begin{subfigure}[b]{0.45\textwidth}
         \centering
         \includegraphics[trim={0 150 0 40},clip,width=\textwidth]{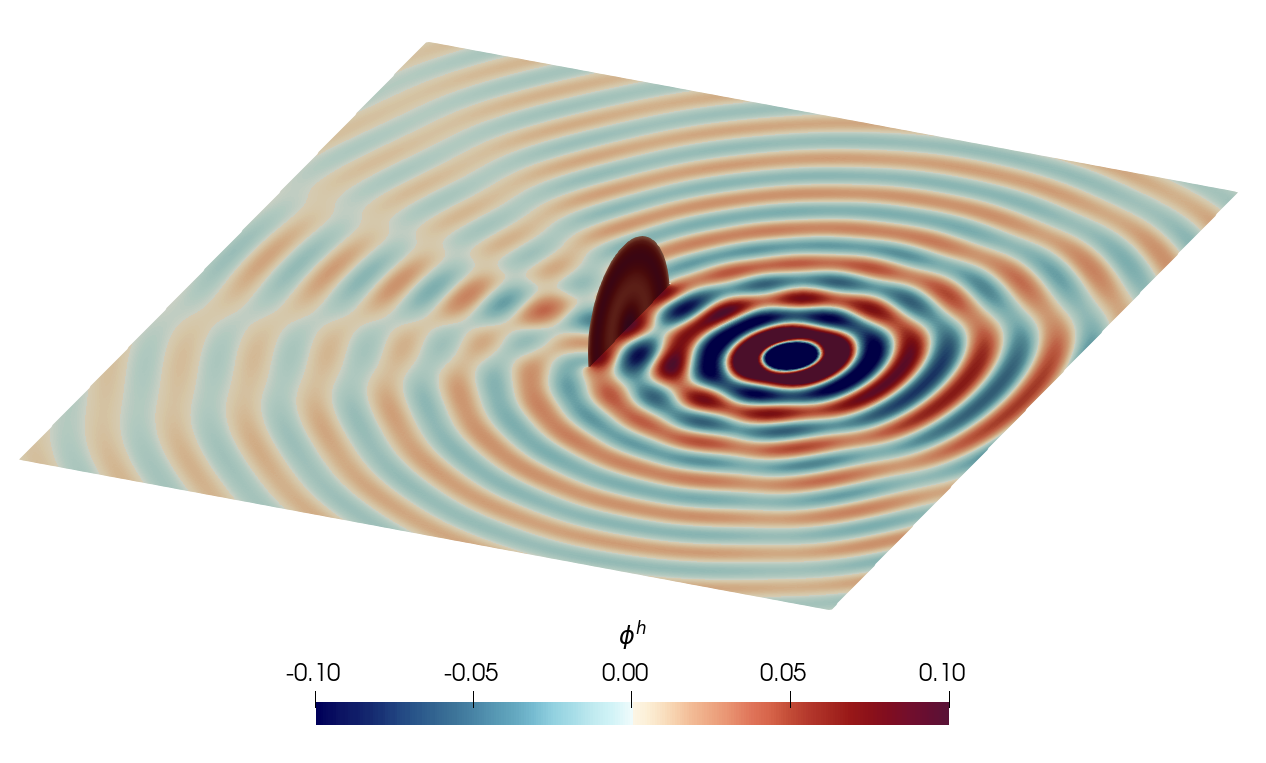}
         \caption{$ka=8$}
         \label{fig:disk_harmonic_snapshot_k8}
     \end{subfigure}
    \begin{subfigure}[b]{0.45\textwidth}
         \centering
         \includegraphics[trim={0 150 0 40},clip,width=\textwidth]{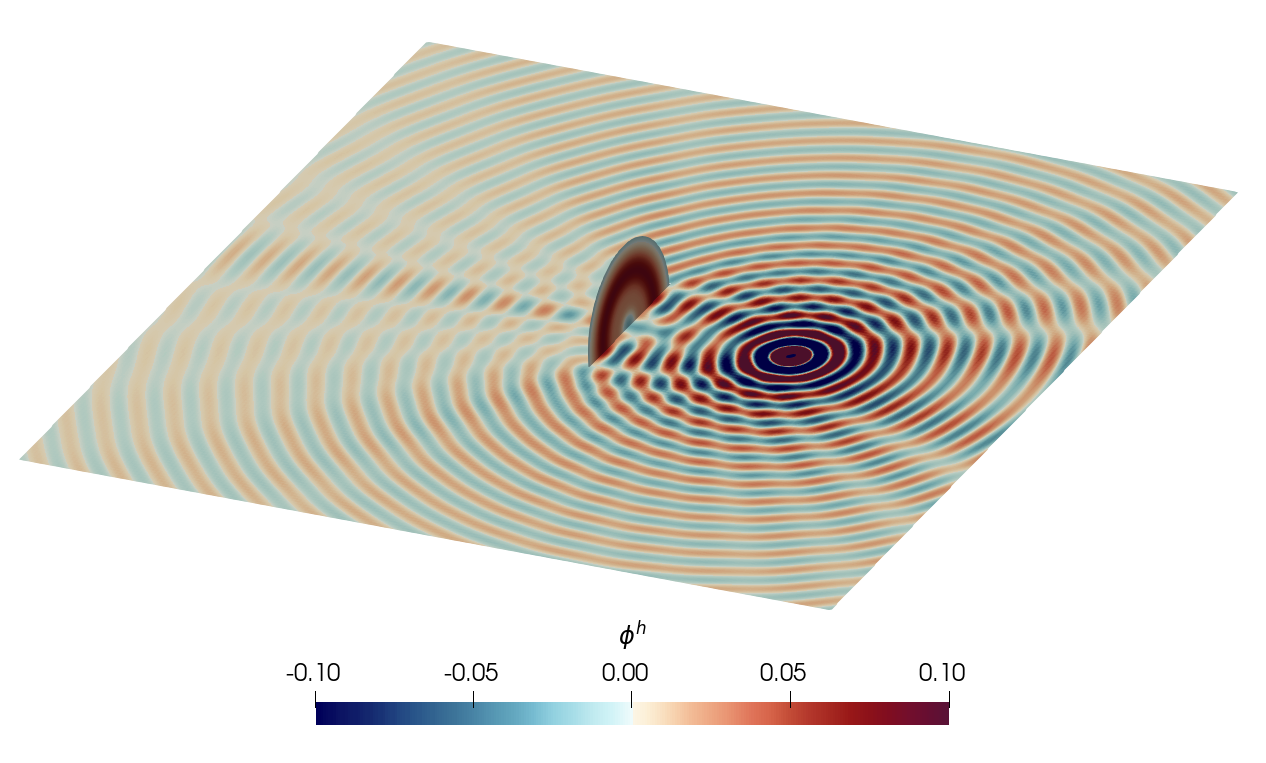}
         \caption{$ka=16$}
         \label{fig:disk_harmonic_snapshot_k16}
     \end{subfigure}
    \caption{Visualizations of the predicted acoustic field of a harmonic point source scattered by a disk at selected frequencies: (a) $ka=2$; (b) $ka=4$; (c) $ka=8$; (d) $ka=16$. The instantaneous acoustic potential is shown at nondimensional time $t=50$.}
    \label{fig:disk_harmonic_snapshots}
\end{figure}

\begin{figure}[h]
    \centering
    \begin{subfigure}[b]{0.45\textwidth}
         \centering
         \includegraphics[trim={0 3 0 3},clip,width=\textwidth]{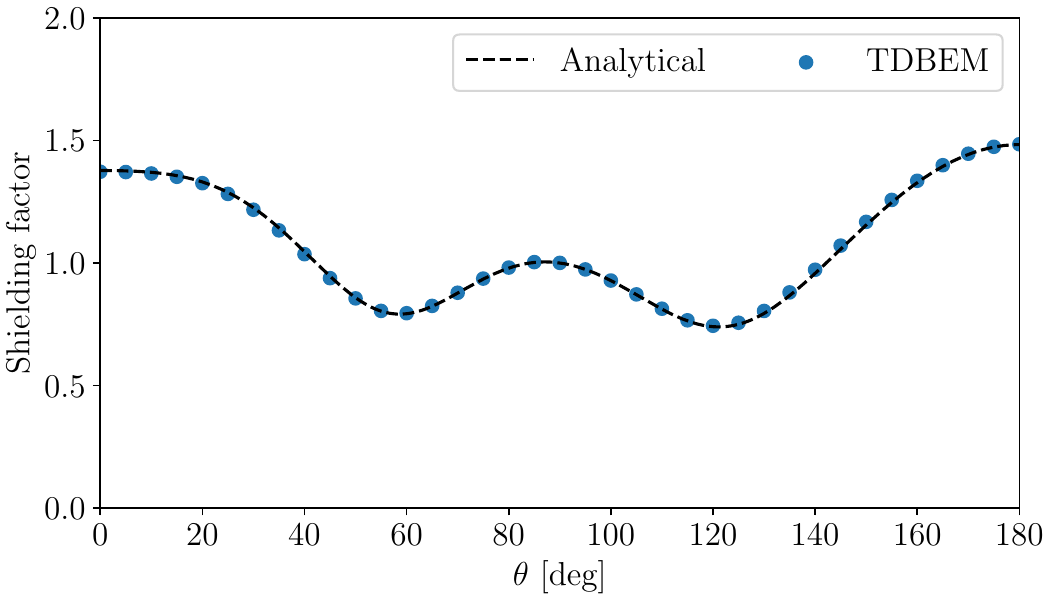}
         \caption{$ka=2$}
         \label{fig:disk_harmonic_directivity_k2}
     \end{subfigure}
    \begin{subfigure}[b]{0.45\textwidth}
         \centering
         \includegraphics[trim={0 3 0 3},clip,width=\textwidth]{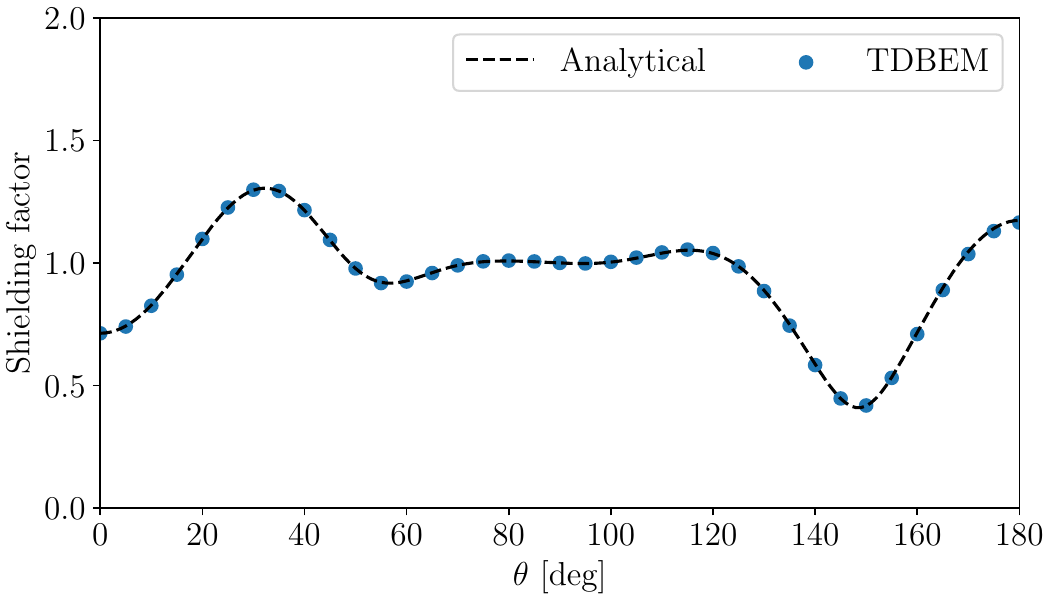}
         \caption{$ka=4$}
         \label{fig:disk_harmonic_directivity_k4}
     \end{subfigure}\\
    \begin{subfigure}[b]{0.45\textwidth}
         \centering
         \includegraphics[trim={0 3 0 3},clip,width=\textwidth]{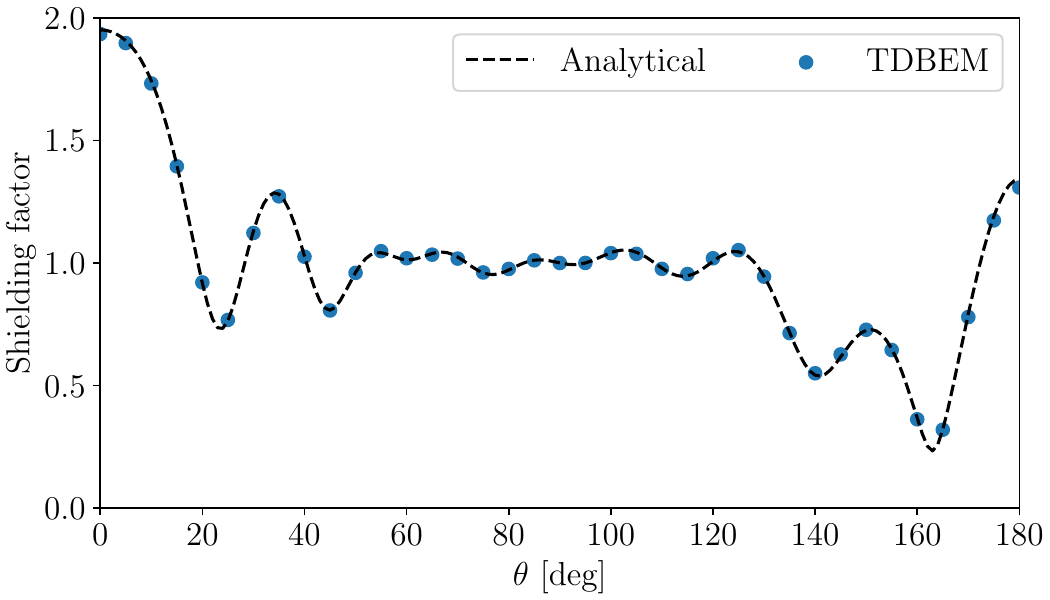}
         \caption{$ka=8$}
         \label{fig:disk_harmonic_directivity_k8}
     \end{subfigure}
    \begin{subfigure}[b]{0.45\textwidth}
         \centering
         \includegraphics[trim={0 3 0 3},clip,width=\textwidth]{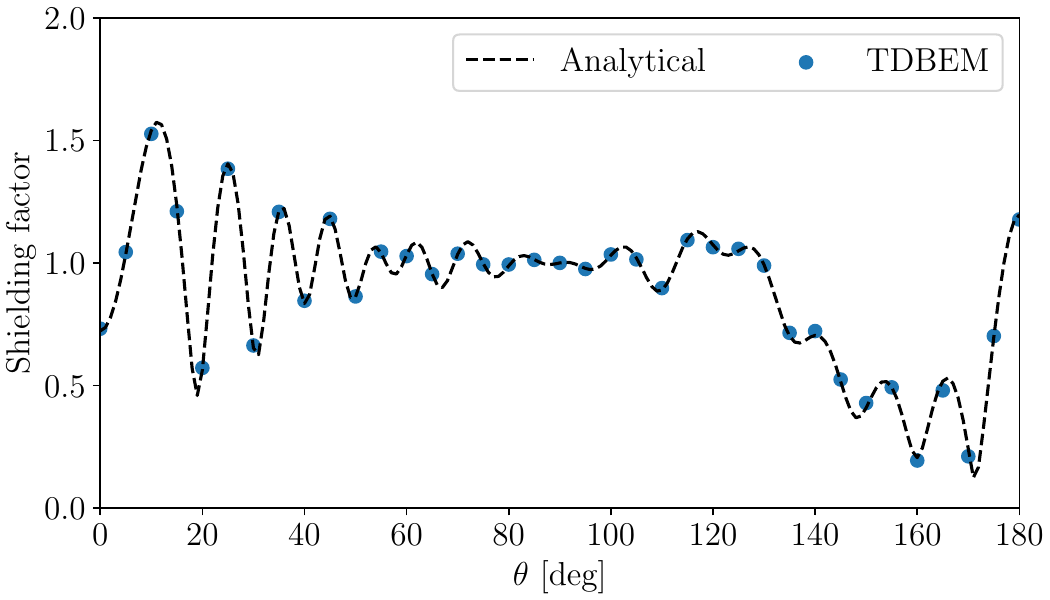}
         \caption{$ka=16$}
         \label{fig:disk_harmonic_directivity_k16}
     \end{subfigure}
    \caption{Comparison of predicted and analytical shielding factors for scattering of a harmonic point source by a disk at selected frequencies: (a) $ka=2$; (b) $ka=4$; (c) $ka=8$; (d) $ka=16$.}
    \label{fig:disk_harmonic}
\end{figure}

For this case, we focus on a constant harmonic frequency forcing ($Q(t) = -\cos(kt)$), and compare with the analytical solution given by Bowman et al.~\cite{bowman1969electromagnetic}. As with the sphere case, the numerical data are collected after the start of the simulation when initial transient effects have exited the relevant parts of the domain. The acoustic potential distributions $\phi^h$ for harmonic forcing at $ka = 2,4,8,16$ are shown in Fig. \ref{fig:disk_harmonic_snapshots} on the disk surface and on an observer plane $x \in [-5,5]$, $y \in [-5,5]$. Compared to the sphere case, the visualizations show substantial qualitative differences with a stronger influence of the scattering body along the $x$-axis due to the flat geometry.

The shielding factor directivities on the circular observer array are compared with the analytical solution for $ka = 2,4,8,16$ in Fig. \ref{fig:disk_harmonic}. As with the sphere case, the agreement between the predicted and analytical shielding factors is strong. The solutions are essentially coincident at lower frequencies at all polar angles, while small numerical errors appear at higher frequencies due to the fixed resolution.

\FloatBarrier

\subsection{Plane}
\begin{figure}[h]
\centering
\begin{tikzpicture}

\def\angEl{28} 
\def\angAz{-112} 

\tikzset{xyplane/.style={cm={cos(\angAz),sin(\angAz)*sin(\angEl),-sin(\angAz),cos(\angAz)*sin(\angEl),(0,0)}}}

\draw[xyplane,fill=gray!40] (-2.5,-2.5) -- (-2.5,2.5) -- (2.5,2.5) -- (2.5,-2.5) -- cycle;

\draw[xyplane,->] (-2,-5) -- (-2,-4) node[right] {$M_\infty$};

\draw[xyplane,<-] (0,-4) node[left] {$-x$} -- (0,0);
\draw[xyplane,->] (0,0) -- (0,4) node[right] {$x$};

\draw[xyplane,<-] (-4,0) node[above right] {$y$} -- (0,0);
\draw[xyplane,->] (0,0) -- (4,0) node[below left] {$-y$};

\draw[xyplane] (2.5,-2.5) node[left] {($-1,-1$)};
\draw[xyplane] (-2.5,-2.5) node[above] {($-1,1$)};
\draw[xyplane] (2.5,2.5) node[below] {($1,-1$)};
\draw[xyplane] (-2.5,2.5) node[right] {($1,1$)};

\draw[->] (0,0) -- (0,2.2) node[above] {$z$};

\draw (0,0.4) -- (1,0.4);
\draw (0,0) -- (1,0);
\draw[<->,dashed] (0.9,0) -- (0.9,0.4) node[pos=0.5,right,inner sep=2pt] {$0.2$};
\fill[blue] (0,0.4) circle (2pt) node[left,inner sep=2pt] {$R_{src}$};

\draw (0,1) -- (0.8,1);
\draw[<->,dashed] (0.7,0) -- (0.7,1) node[pos=0.7,right,inner sep=2pt] {$0.5$};
\fill[red] (0,1) circle (2pt) node[left,inner sep=2pt] {$R_{obs}$};
\end{tikzpicture}
\caption{Plane geometry, showing the source (blue dot, $R_{src}$) and observer (red dot, $R_{obs}$) locations, with the uniform background flow in the positive $x$ direction.}
\label{fig:plate_geo}
\end{figure}
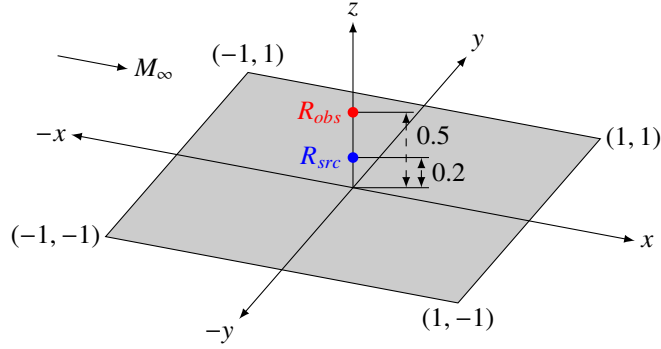

The final validation case is acoustic scattering of a transient source in a uniform flow parallel to a flat plane. The analytical solution for scattering by an infinite plane is easily constructed by the method of images. The infinite plane problem can be approximated numerically for short times with a sufficiently large finite rectangular plate geometry, as illustrated in Fig. \ref{fig:plate_geo}. Although the finite plate will introduce edge diffraction effects that are not present in the infinite plane case, causality ensures these effects can only appear in the solution after the incident wave propagates from the source point to an edge and subsequently to the observer point. Therefore, we can choose a short time interval in which the numerical case with a finite geometry is directly comparable to the analytical solution for an infinite geometry.

We define a transient source strength function:
\begin{equation}
    Q(t) = \frac{t-t_0}{\sigma}\exp\left(-\frac{(t-t_0)^2}{2\sigma^2} + \frac{1}{2}\right).
\end{equation}
This function is proportional to the derivative of a Gaussian function and is normalized to unit amplitude. For this case, we choose $t_0 = 2$ and $\sigma = 0.085$.

The mean background flow is uniform and parallel to the positive $x$-axis. We therefore apply the PGL transformation in the TDBEM computation to account for background flow effects:
\begin{subequations}
\begin{align}
\tilde{x} & = \frac{x}{\beta} + \frac{M_\infty \cdot x}{\beta^2(1+\beta)}M_\infty, \\
\tilde{t} & = t + \frac{M_\infty \cdot x}{\beta^2}.
\end{align}
\end{subequations}

\begin{figure}[h]
    \centering
    \includegraphics[trim={0 0 0 615},clip,width=0.6\textwidth]{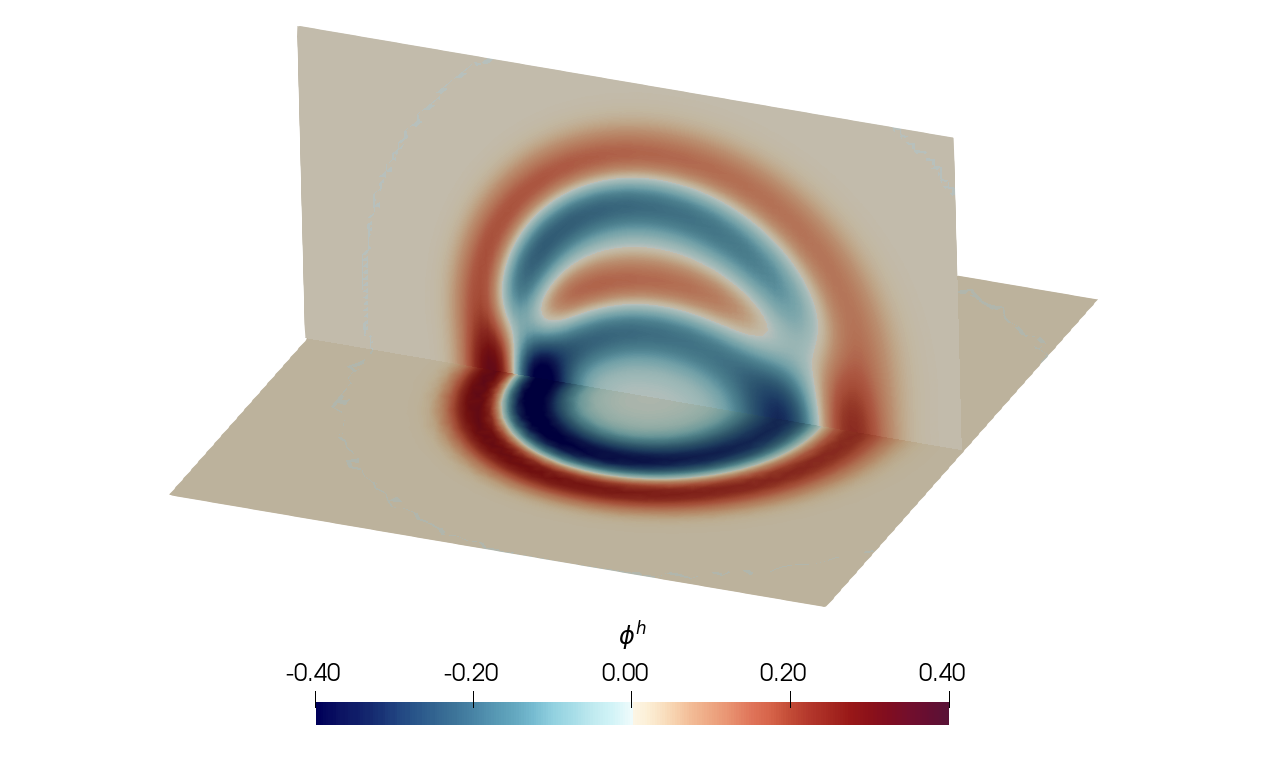}\\
    \begin{subfigure}[b]{0.45\textwidth}
         \centering
         \includegraphics[trim={100 155 100 20},clip,width=\textwidth]{plate_M02_250.png}
         \caption{$|M_\infty| = 0.2$}
         \label{fig:plate_snapshot_M02}
     \end{subfigure}
    \begin{subfigure}[b]{0.45\textwidth}
         \centering
         \includegraphics[trim={100 155 100 20},clip,width=\textwidth]{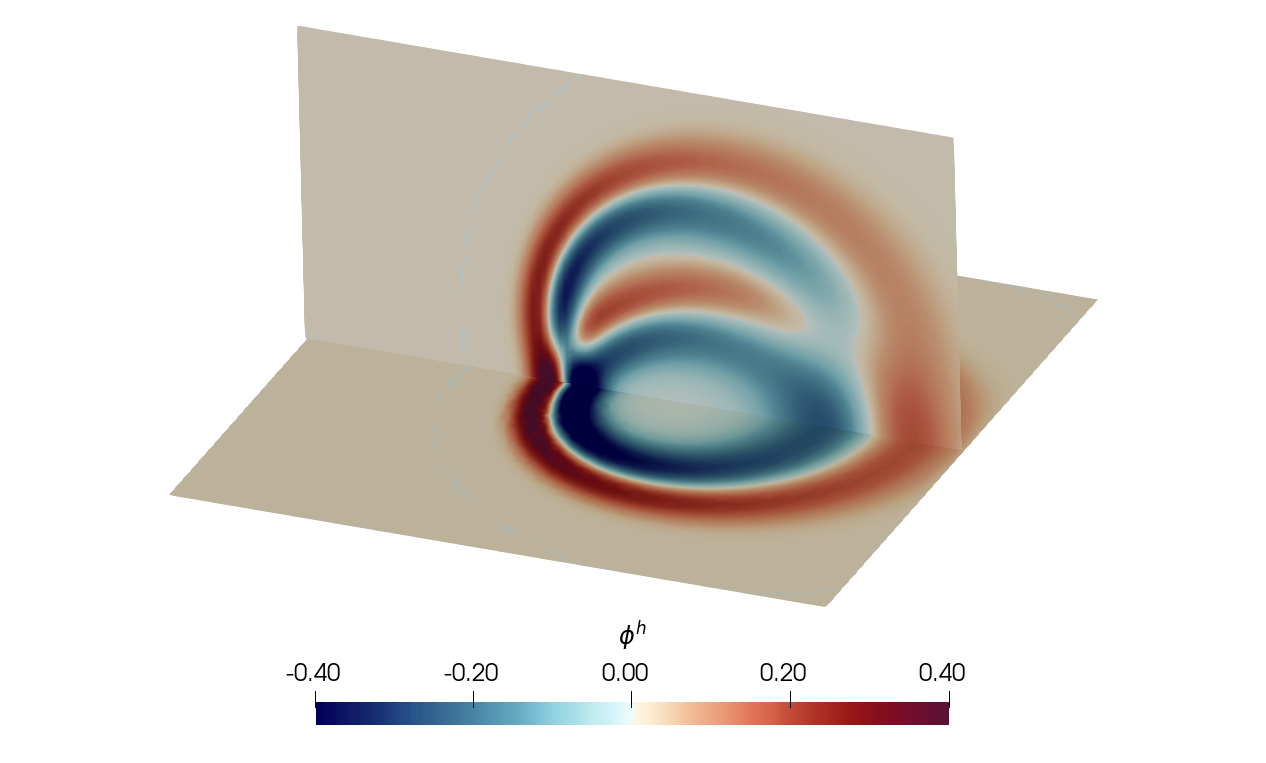}
         \caption{$|M_\infty| = 0.5$}
         \label{fig:plate_snapshot_M05}
     \end{subfigure}
    \caption{Visualizations of the predicted acoustic potential field of a transient point source scattered by a plane with a uniform mean flow: (a) $|M_\infty| = 0.2$; (b) $|M_\infty| = 0.5$. The instantaneous acoustic potential is shown at nondimensional time $t=2.5$.}
    \label{fig:plate_snapshots}
\end{figure}

\begin{figure}[h]
    \centering
    \begin{subfigure}[b]{0.45\textwidth}
         \centering
         \includegraphics[trim={0 5 0 2},clip,width=\textwidth]{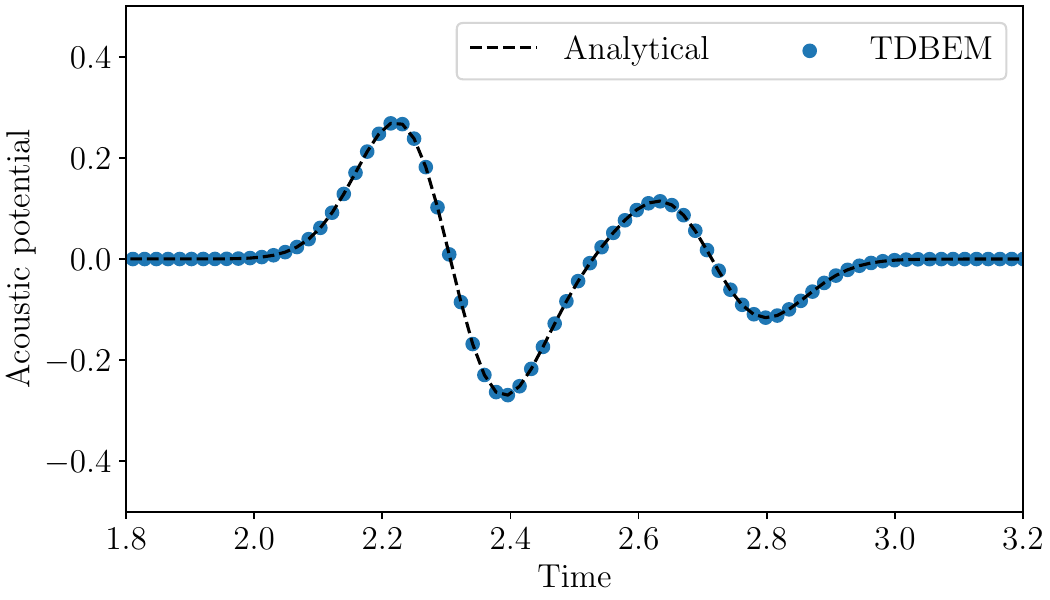}
         \caption{$|M_\infty| = 0.2$}
         \label{fig:plate_signal_M02}
     \end{subfigure}
    \begin{subfigure}[b]{0.45\textwidth}
         \centering
         \includegraphics[trim={0 5 0 2},clip,width=\textwidth]{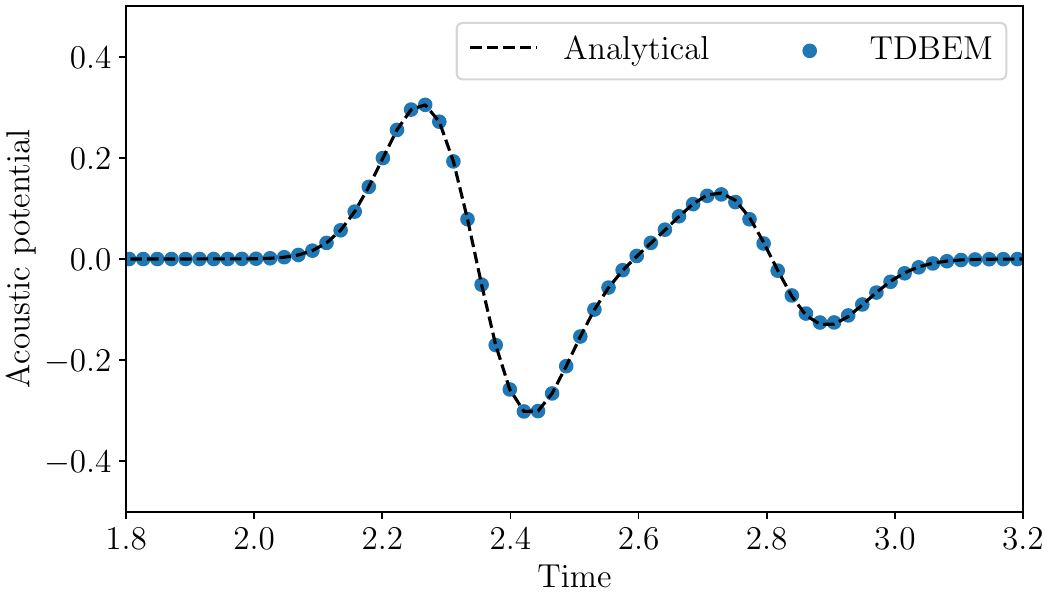}
         \caption{$|M_\infty| = 0.5$}
         \label{fig:plate_signal_M05}
     \end{subfigure}
    \caption{Comparison of the predicted and analytical acoustic potential due to a transient point source scattered by a plane with a uniform mean flow: (a) $M_\infty = (0.2, 0, 0)$; (b) $M_\infty = (0.5, 0, 0)$. The acoustic potential signal is shown at the observer point $R_{obs} = (0,0,0.5)$.}
    \label{fig:plate_signals}
\end{figure}

The finite plate is discretized using a single layer of elements with $6282$ faces and $3242$ nodes, and computations are performed at $\mbox{CFL} = 0.5$. Note that the PGL transform scales the spatial coordinates, so specifying a constant $\mbox{CFL}$ value results in slightly different timestep sizes at different Mach numbers for the same discretization in physical space. The instantaneous acoustic potential distribution is visualized at $t = 2.5$ for Mach numbers $|M_\infty| = 0.2$ and $|M_\infty| = 0.5$ in Fig. \ref{fig:plate_snapshots} on the plane surface and on a central observer plane aligned with the flow. The upstream amplification due to the background flow is clearly visible as the Mach number increases. The predicted and analytical acoustic potential signals at $R_{obs} = (0,0,0.5)$ are compared for both cases in Fig. \ref{fig:plate_signals}. The signals are essentially coincident for all timesteps at both Mach numbers, demonstrating strong agreement.

\FloatBarrier

\section{Application to a trailing edge-mounted propeller}
\begin{figure}[h]
\centering
\begin{tikzpicture}[scale=3.5]
\draw[ultra thick] (-0.2,-0.2) -- (0.5,-0.2);

\fill[black] (0,0) circle (0.5pt);
\draw[thick, fill=gray!40] (0,0.05) ellipse (0.01 and 0.05);
\draw[thick, fill=gray!40] (0,-0.05) ellipse (0.01 and 0.05);

\draw (0.025,-0.1) -- (0.175,-0.1);
\draw (0.025,0.1) -- (0.175,0.1);
\draw[<->,dashed] (0.15,-0.1) -- (0.15,0.1) node[pos=0.5,right] {$D$};

\draw (-0.025,0) -- (-0.275,0);
\draw (-0.225,-0.2) -- (-0.275,-0.2);
\draw[<->,dashed] (-0.25,-0.2) -- (-0.25,0) node[pos=0.5,left] {$H$};

\draw (0,-0.225) -- (0,-0.275);
\draw (-0.2,-0.225) -- (-0.2,-0.275);
\draw[<->,dashed] (0,-0.25) -- (-0.2,-0.25) node[pos=0.5,below] {$L$};

\draw[thick] (0.94*0.75,0.34*0.75) arc (20:150:0.75)
node[pos=0.55,above, inner sep=2pt] {Reflection side array}
node[pos=0,below, inner sep=2pt] {$20^\circ$}
node[pos=1,below, inner sep=2pt] {$150^\circ$};
\draw[thick] (0.94*0.75,-0.34*0.75) arc (-20:-150:0.75)
node[pos=0.55,below, inner sep=2pt] {Shadow side array}
node[pos=0,above, inner sep=2pt] {$340^\circ$}
node[pos=1,above, inner sep=2pt] {$210^\circ$};

\draw (0.025,0) -- (0.125,0);
\draw (0.25,0) -- (0.375,0);
\draw[<->,dashed] (0.35,0) arc (0:60:0.35) node[pos=0.5,right,inner sep=2pt] {$\theta$};

\draw[<->,dashed] (0,0) -- (0.5*0.75,0.87*0.75) node[pos=0.7,left,inner sep=2pt] {$r = 1.75$ m};

\end{tikzpicture}
\caption{Schematic of trailing edge mounted propeller configuration, showing the position of the plate (thick line), and reflection side and shadow side microphone arrays relative to the propeller.}
\label{fig:prop_geo}
\end{figure}
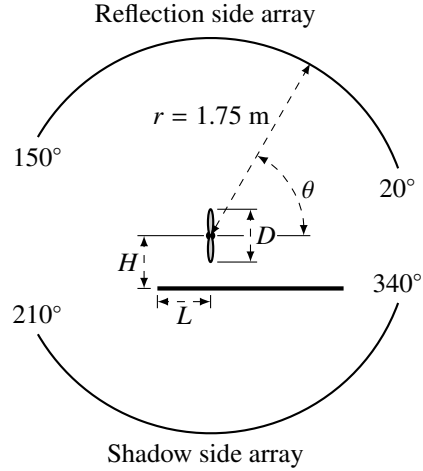

We now apply our scattering methodology to a propeller installed near the trailing edge of a flat plate. This case replicates an experiment performed by Hanson et al.~\cite{Hanson2022}, allowing us to compare our numerical predictions with experimental measurements. The trailing edge-mounted propeller geometry is illustrated in Fig. \ref{fig:prop_geo}. A single two-bladed propeller is mounted near the trailing edge of a rectangular flat plate. The propeller has a pitch and diameter of $D = 0.2286$~m and is operated in a hover condition at a rotation rate of $\Omega = 7000$~RPM. Two rectangular plates are considered: a long plate with a chord of $1.7$~m, and with a short plate with a chord of $0.61$~m, each with a span of $1.5$~m and a thickness of $0.012$~m. Each plate is placed at a range of normal and chordwise positions from $H/D = 0.6$ to $H/D = 1$ and $L/D = -1$ to $L/D = 2$, all at mid-span, where $H$ is the normal separation from the propeller axis to the plate and $L$ is the chordwise distance from the propeller disk plane to the trailing edge. Farfield noise measurements were collected on two 14-microphone arcs centered on the propeller hub at a radius of $1.75$~m, spanning polar angles from $20^\circ$ to $150^\circ$ on the reflection side and $210^\circ$ to $340^\circ$ on the shadow side of the plate.

The incident acoustic field is obtained from a Reynolds-averaged Navier-Stokes (RANS) CFD simulation with the Ffowcs Williams-Hawkings (FWH) acoustic analogy~\cite{ffowcs1969sound}. The experimental results show minimal aerodynamic interaction between the propeller and the plate in all operating conditions, and the measured installation effects are attributed to acoustic scattering and shielding~\cite{Hanson2022}. Therefore, a single simulation of the isolated propeller is used to provide the incident field for all cases. The flow solution is computed in a rotating reference frame on a very fine 251 million grid point mesh using the open-source SU2 solver~\cite{economon2016su2}. The acoustic field is evaluated at 2$^\circ$ blade azimuth per timestep.

The $1.7$~m plate is discretized with $5398$ faces and $2791$ nodes and the $0.61$~m plate is discretized with $1938$ faces and $1030$ nodes, resulting in $\mbox{CFL} \approx 0.5$ for both geometries. Because the experiment was conducted in a hover condition, the background flow is assumed to be at rest and no variable transformation is applied. The simulations were run for $H/D = 0.6, 0.75, 1$ and $L/D = -1, 0, 1, 2$. To compare against the experimental results with significant broadband content not captured by the RANS solution, we focus on the sound pressure level (SPL) at the blade passing frequency (BPF).

\begin{figure}[h]
    \centering
    \begin{subfigure}[b]{0.45\textwidth}
         \centering
         \includegraphics[width=\textwidth]{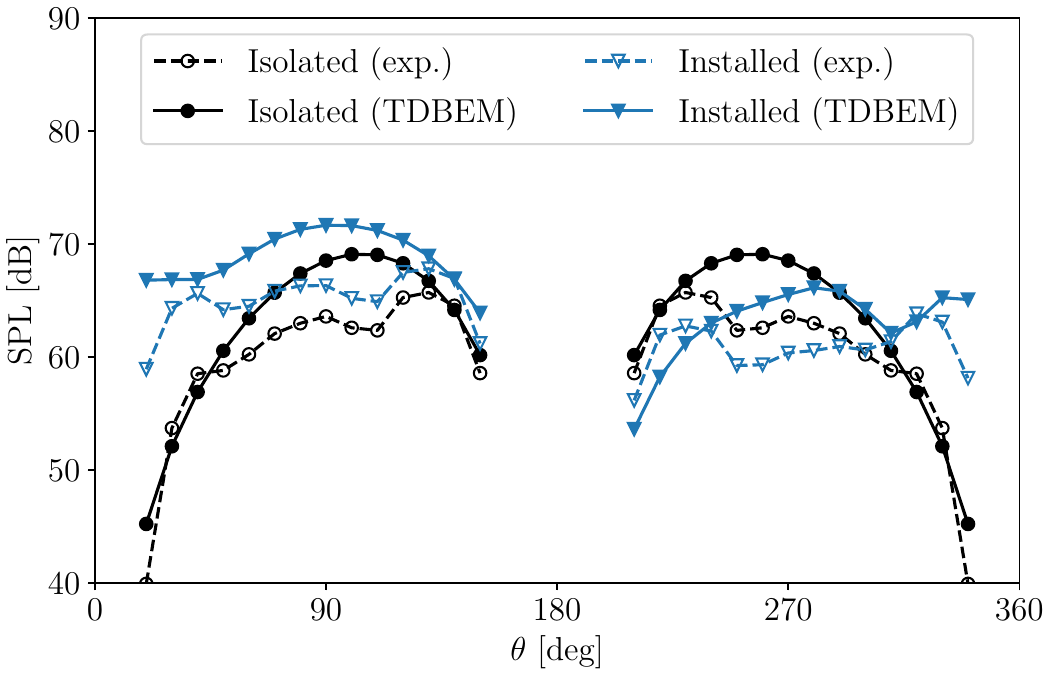}
         \caption{Long plate, $H/D=1$, $L/D=0$}
     \end{subfigure}
    \begin{subfigure}[b]{0.45\textwidth}
         \centering
         \includegraphics[width=\textwidth]{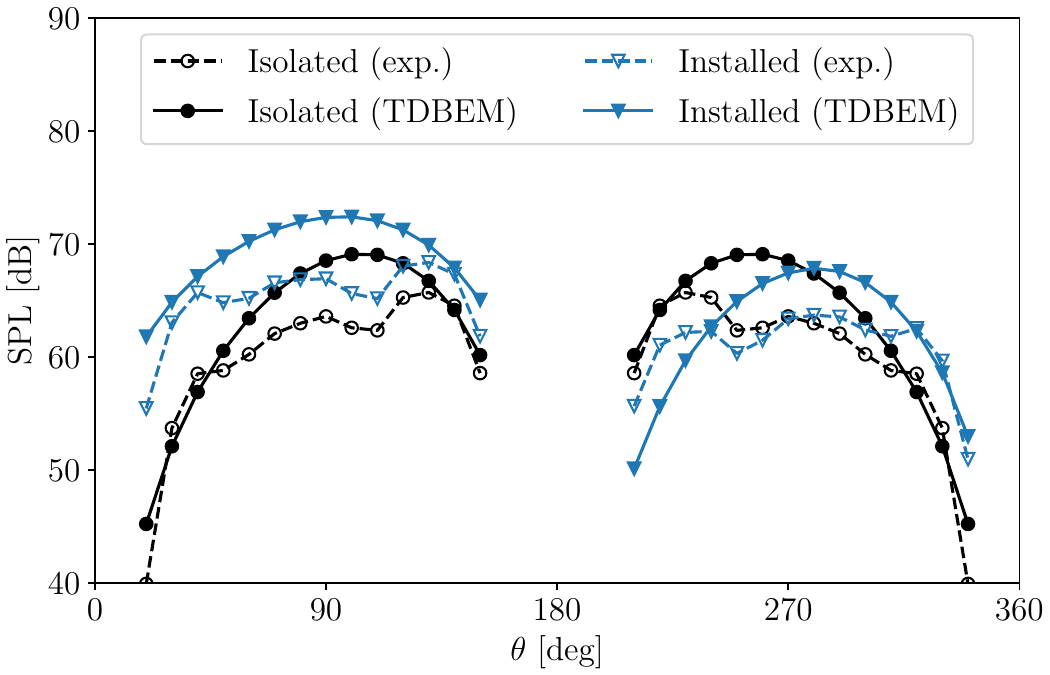}
         \caption{Short plate, $H/D=1$, $L/D=0$}
     \end{subfigure}\\
    \begin{subfigure}[b]{0.45\textwidth}
         \centering
         \includegraphics[width=\textwidth]{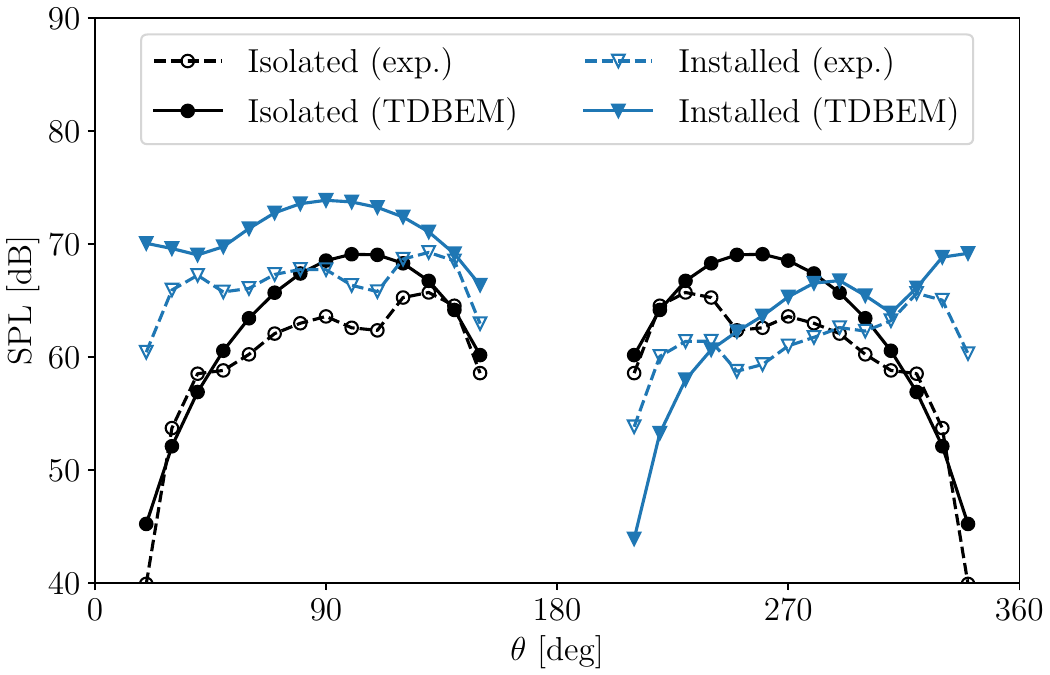}
         \caption{Long plate, $H/D=0.75$, $L/D=0$}
     \end{subfigure}
    \begin{subfigure}[b]{0.45\textwidth}
         \centering
         \includegraphics[width=\textwidth]{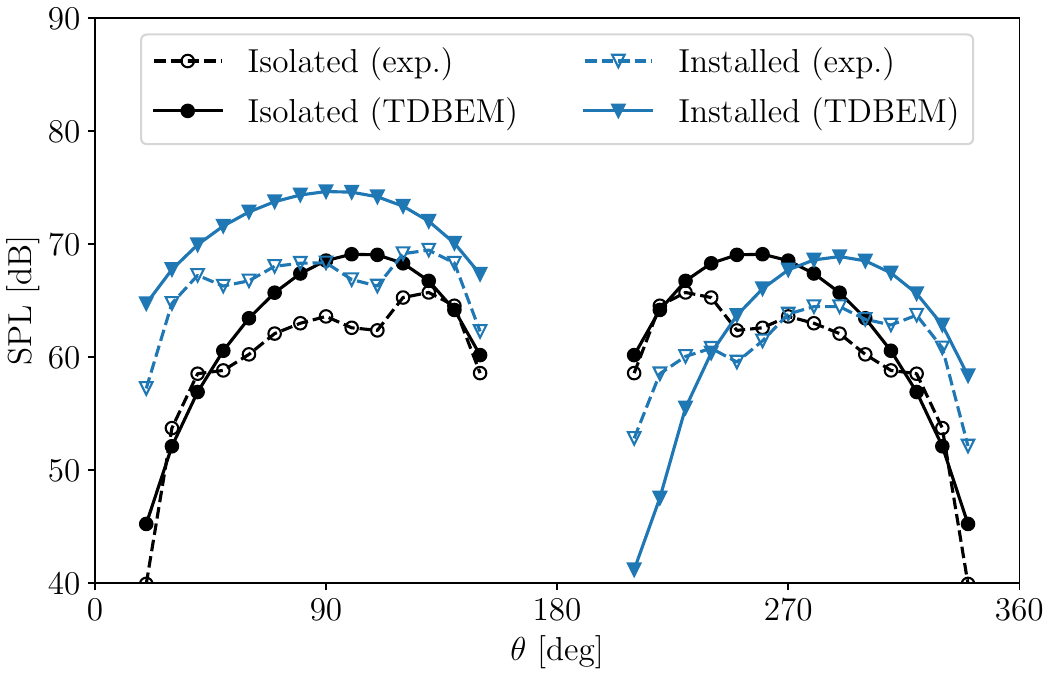}
         \caption{Short plate, $H/D=0.75$, $L/D=0$}
     \end{subfigure}\\
    \begin{subfigure}[b]{0.45\textwidth}
         \centering
         \includegraphics[width=\textwidth]{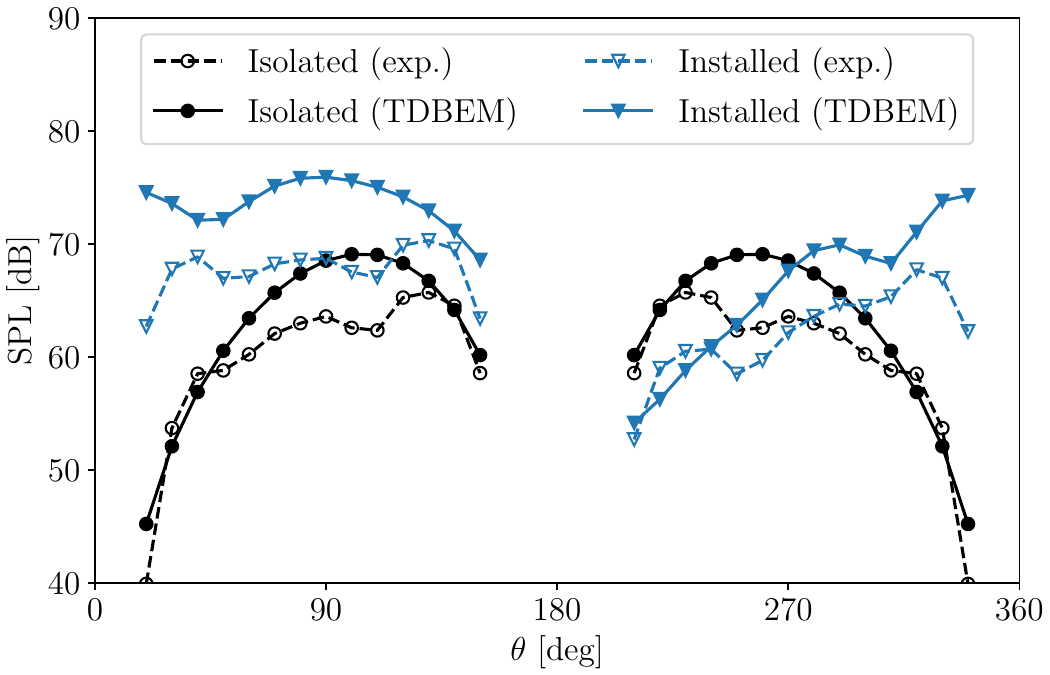}
         \caption{Long plate, $H/D=0.6$, $L/D=0$}
     \end{subfigure}
    \begin{subfigure}[b]{0.45\textwidth}
         \centering
         \includegraphics[width=\textwidth]{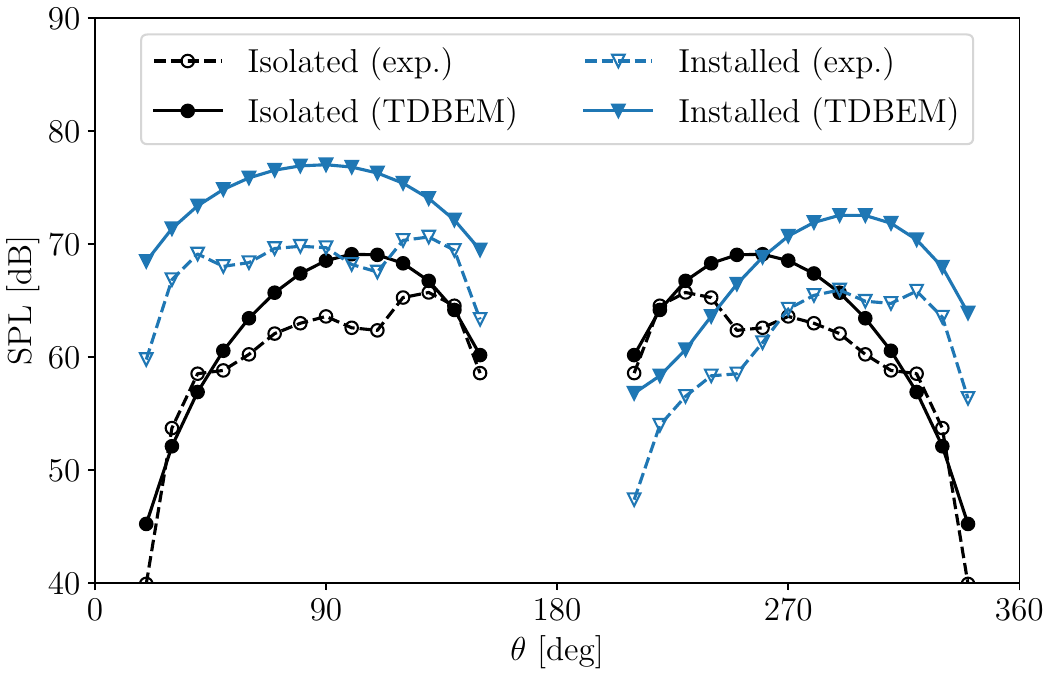}
         \caption{Short plate, $H/D=0.6$, $L/D=0$}
     \end{subfigure}
    \caption{Comparison of predicted (solid markers and lines) and measured (open markers and dashed lines) directivites for SPL at BPF for $L/D = 0$ (propeller disk aligned with trailing edge). (a) long plate, $H/D=1$; (b) short plate, $H/D=1$; (c) long plate, $H/D=0.75$; (d) short plate, $H/D=0.75$; (e) long plate, $H/D=0.6$; (f) short plate, $H/D=0.6$.}
    \label{fig:prop_SPL_directivity}
\end{figure}

\begin{figure}[h]
    \centering
    \begin{subfigure}[b]{0.45\textwidth}
         \centering
         \includegraphics[width=\textwidth]{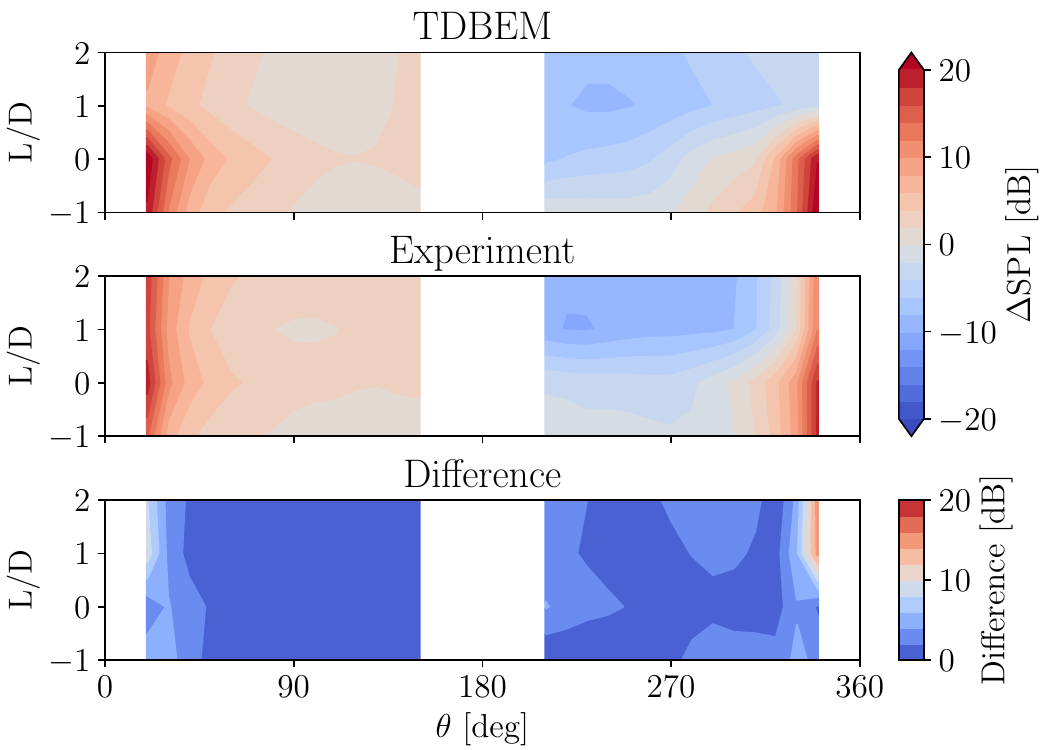}
         \caption{Long plate, $H/D=1$}
     \end{subfigure}
    \begin{subfigure}[b]{0.45\textwidth}
         \centering
         \includegraphics[width=\textwidth]{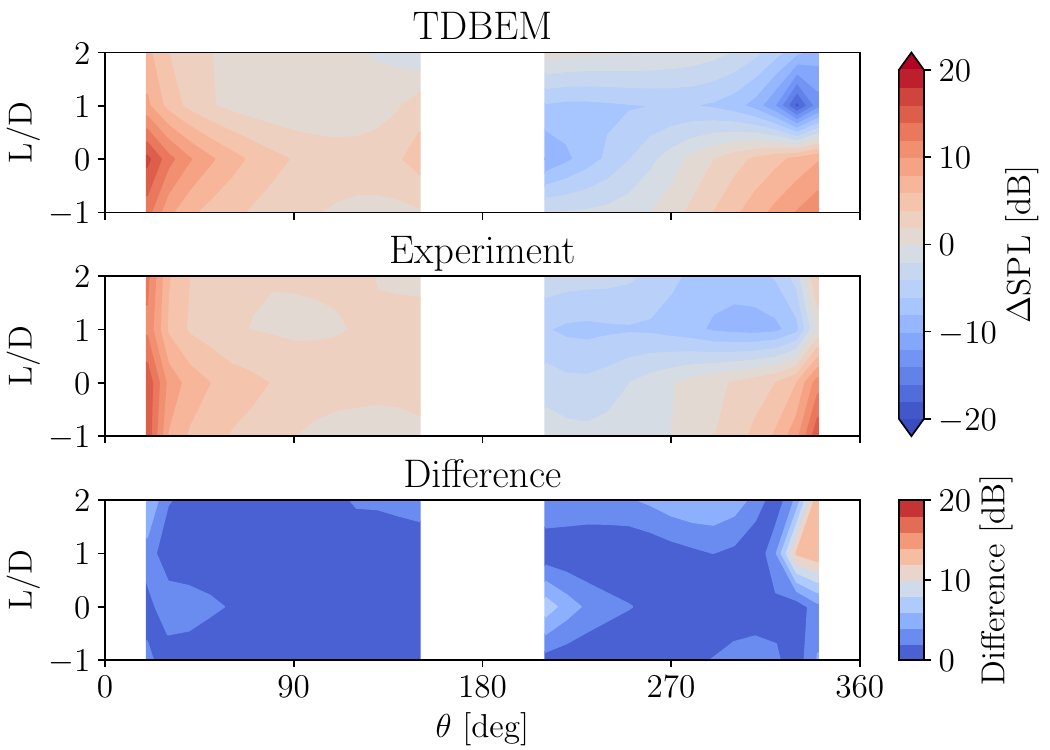}
         \caption{Short plate, $H/D=1$}
     \end{subfigure}\\
    \begin{subfigure}[b]{0.45\textwidth}
         \centering
         \includegraphics[width=\textwidth]{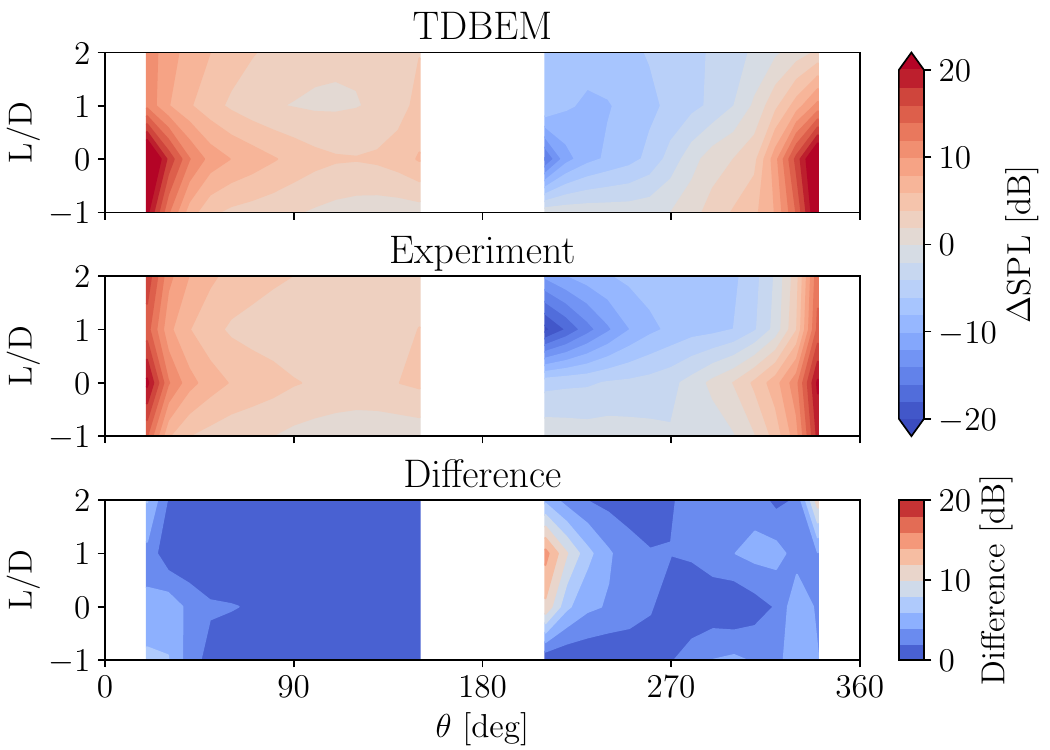}
         \caption{Long plate, $H/D=0.75$}
     \end{subfigure}
    \begin{subfigure}[b]{0.45\textwidth}
         \centering
         \includegraphics[width=\textwidth]{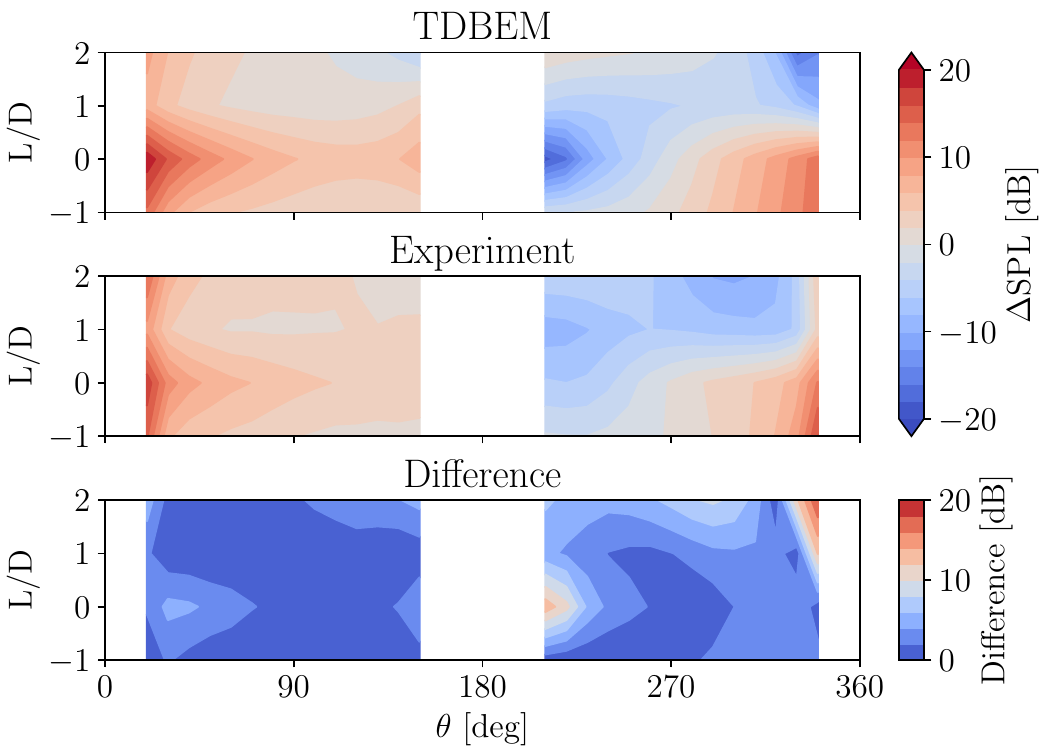}
         \caption{Short plate, $H/D=0.75$}
     \end{subfigure}\\
    \begin{subfigure}[b]{0.45\textwidth}
         \centering
         \includegraphics[width=\textwidth]{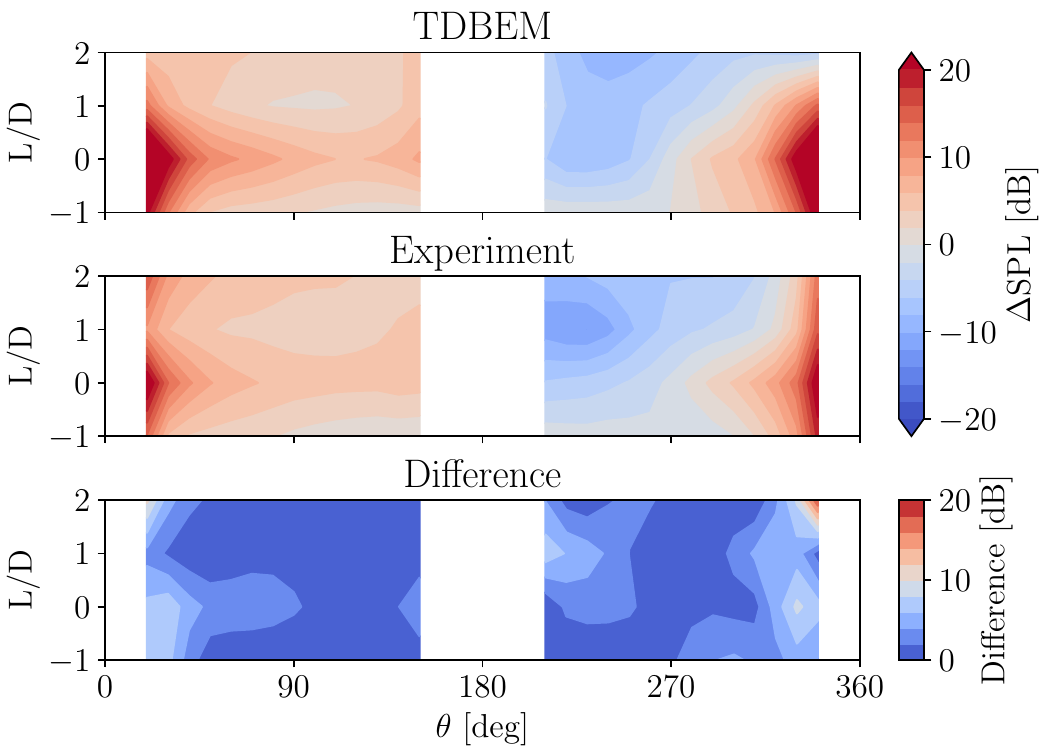}
         \caption{Long plate, $H/D=0.6$}
     \end{subfigure}
    \begin{subfigure}[b]{0.45\textwidth}
         \centering
         \includegraphics[width=\textwidth]{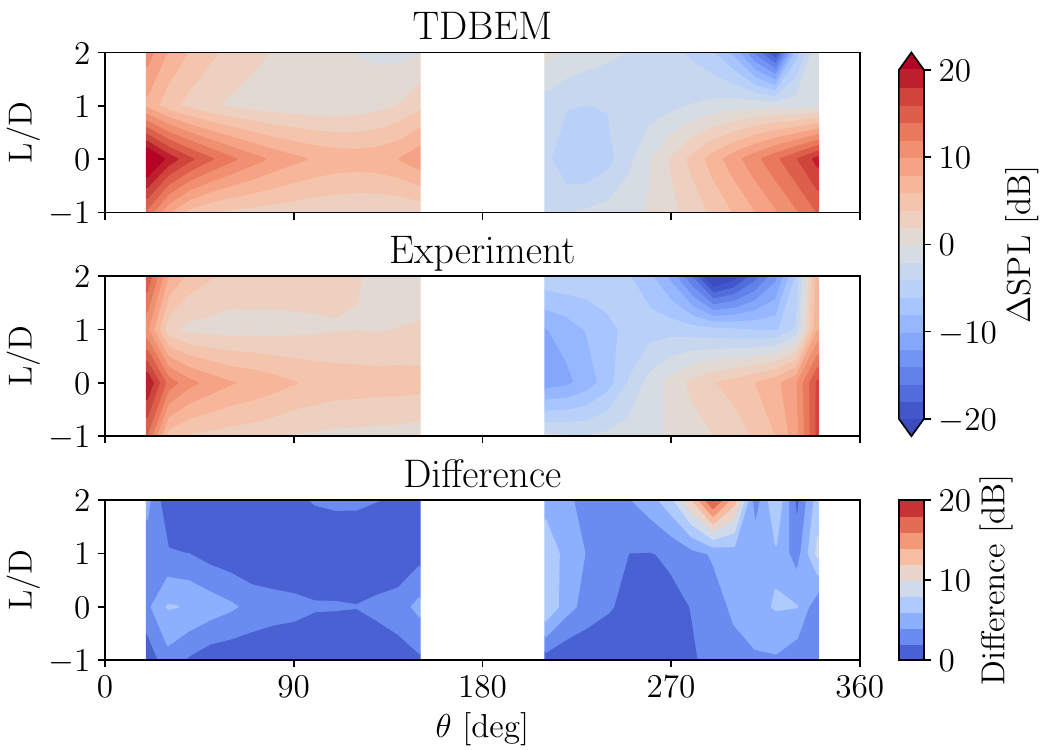}
         \caption{Short plate, $H/D=0.6$}
     \end{subfigure}
    \caption{Comparison of $\Delta$SPL between the installed and isolated configurations with varying $L/D$. (a) long plate, $H/D=1$; (b) short plate, $H/D=1$; (c) long plate, $H/D=0.75$; (d) short plate, $H/D=0.75$; (e) long plate, $H/D=0.6$; (f) short plate, $H/D=0.6$. Each subfigure shows contours of the predicted $\Delta$SPL (top), measured $\Delta$SPL (middle), and the magnitude of the difference between the predicted and measured $\Delta$SPL (bottom).}
    \label{fig:prop_delta_SPL}
\end{figure}

The predicted and measured SPL directivities at the BPF are compared in Fig. \ref{fig:prop_SPL_directivity} for $L/D = 0$, which corresponds to the propeller disk aligned with the trailing edge. In each plot, the isolated directivities are included alongside the installed directivities to highlight the acoustic installation effect. We note some systematic differences in the shape of the isolated propeller directivity, with overprediction of the measured values by the RANS-FWH numerical approach for intermediate polar angles around $90^\circ$ and $270^\circ$. These differences mostly carry over from the incident field prediction into the installed propeller predictions, which show very similar overprediction around $90^\circ$ and $270^\circ$. These discrepancies may be due to aerodynamic and acoustic differences between the numerical and experimental case setups, including the spinner, motor, and test stand, which are present in the experiment and excluded from the numerical geometry, as well as the presence of a wind tunnel nozzle and collector in the anechoic chamber, which potentially act as additional scattering surfaces in the experimental setup. They may also be related to separated flow phenomena on the propeller blades that are not sufficiently captured in the RANS flow solution.

Despite differences in the absolute levels due to the incident field prediction, the separation between the isolated and installed levels due to acoustic scattering and shielding is very well captured with the TDBEM. On the reflection side, scattering by the plate results in amplification across the entire observer array, which is most pronounced at the forward end (around $20^\circ$). The amplification is also stronger for smaller $H/D$. The long plate results in more amplification at the forward end of the observer array compared to the short plate, although the numerical predictions show a stronger increase for the long plate between $20^\circ$ and $40^\circ$ than the experimental measurements. On the shadow side, the scattering effect is more complicated, with amplification at forward angles (around $340^\circ$) and reduction at backward angles (around $210^\circ$), and a crossover around $270^\circ$. The crossover location is particularly well captured as it moves forward with increasing $H/D$ and for the long plate configuration. Other key parametric trends are again captured, such as increased amplification and reduction for smaller $H/D$. As with the reflection side, the 
predicted amplification is stronger than the measurements for the long plate at the forward end of the observer arc between $320^\circ$ and $340^\circ$.

To further characterize the scattering and shielding effects at different propeller locations, we directly compare the difference between the SPL at the BPF in the installed and isolated configurations over the $H/D$ and $L/D$ parameter ranges, computed as:
\begin{equation}
    \Delta \mbox{SPL} = \mbox{SPL}_{\mbox{installed}} - \mbox{SPL}_{\mbox{isolated}}.
\end{equation}
The predicted and measured $\Delta$SPL are shown for varying $L/D$ in Fig. \ref{fig:prop_delta_SPL}, along with the magnitude of the difference between them. The amplification and reduction patterns are well predicted for each $H/D$. On the reflection side, the amplification is strongest for $L/D = 0$ for both plate sizes and decreases as the propeller disk is moved forward. This decrease is fastest for the short plate. On both the reflection and shadow sides, the strong amplification at the forward observers is significantly reduced as the propeller disk moves in front of the trailing edge ($L/D=0$). The crossover angle between regions of amplification and reduction on the shadow side also moves forward substantially.

The difference between predicted and measured $\Delta$SPL is mostly below 5~dB. The exceptions are almost entirely at the forward and backward ends of the array. We hypothesize that the larger discrepancies at these observer locations may be due to the proximity of the wind tunnel nozzle and collector (which are not acoustically treated) to the microphones in the experimental setup. There is one additional large discrepancy for the short plate at $H/D=0.6$ and $L/D=2$, where the location of the maximum reduction on the shadow side is shifted forward in the prediction relative to the measurement, resulting in a significant difference in $\Delta$SPL around $290^\circ$. Comparing across all angles and positions, the numerical predictions are generally closer to the measurements for larger $H/D$. We hypothesize that these locations show stronger agreement because the aerodynamic isolation approximation is more accurate as $H/D$ is increased.

\FloatBarrier

\section{Conclusions}
In this work, we present a method for efficient prediction of acoustic scattering effects for complex installed sources, such as rotors and propellers. A space-time Galerkin time-domain boundary element method is proposed to enforce the sound-hard acoustic boundary condition on a scattering surface for an arbitrary incident field. We address difficulties associated with the numerical evaluation of the Galerkin double integrals and causal propagation integrals with an efficient decomposition-based approach. The space-time Galerkin formulation leads to a robust, unconditionally stable method that provides good accuracy on coarse discretizations. Unlike other formulations, our method can represent thin geometries with a single layer of elements, greatly reducing the problem size on this class of problems. Mean background flow effects are incorporated with variable transformations, providing flexibility to apply appropriate background flow assumptions for different problems without adding significant complexity.

We present validation results for the TDBEM on three cases with analytical solutions, demonstrating very good agreement. Collectively, the three validation cases capture a range of different features that may appear in realistic aeroacoustic scattering problems, including smoothly curved scattering surfaces and planar surfaces with sharp edges, scattering surfaces enclosing a volume and thin surfaces represented with a single layer of elements, scattering in a medium at rest and with a mean flow, and scattering of harmonic, broadband, and transient sources. The method is accurate for all cases without any tuned numerical parameters as required in other methods. Refinement studies performed on the sphere case demonstrate second order accuracy in space and time for joint refinement and suggest a choice of $\mbox{CFL} \approx 0.5$ to minimize the total error.

Finally, we present an application of the scattering method to a trailing edge-mounted propeller configuration. This case demonstrates successful coupling of the TDBEM with CFD via an acoustic analogy to prescribe the incident field. The numerical results are compared with experimental measurements from Hanson et al.~\cite{Hanson2022}, showing good qualitative and quantitative agreement. In particular, the numerical results accurately capture the amplification and reduction trends due to acoustic scattering over varying plate size and location.

\section*{Acknowledgments}
The authors gratefully acknowledge Filipi Kunz and Shaun Pullin for sharing SU2 CFD simulation data used for the trailing edge-mounted propeller case. The authors also thank Dr. Leonard Lopes and Prof. Qiqi Wang for their valuable suggestions to this work over the period of the research. This research was supported in part by computational resources provided by the Partnership for an Advanced Computing Environment (PACE) at the Georgia Institute of Technology, Atlanta, Georgia, USA (RRID:SCR\_027619).

\appendix
\section{Exact expressions for $I_{\alpha,\beta}$}\label{app:integrals}
The proposed quadrature strategy for the evaluation of the Galerkin integrals requires exact expressions for the following integrals:
\begin{align}
    I_{\alpha,1}(r,z) & = \int_0^r \left(\rho^2 + z^2\right)^{\frac{\alpha}{2}} \rho \, d\rho, \label{eqn:app_I1} \\
    I_{\alpha,2}(r,z) & = \int_0^r \left(\rho^2 + z^2\right)^{\frac{\alpha}{2}} \rho^2 \, d\rho. \label{eqn:app_I2}
\end{align}
Note that $I_{\alpha,\beta}$ for $\beta>2$ can be expressed as a combination of Eqs. (\ref{eqn:app_I1}) and (\ref{eqn:app_I2}).

For Eq. (\ref{eqn:app_I1}), let $R = \sqrt{\rho^2 + z^2}$. Then:
\begin{equation}
    I_{\alpha,1}(r,z) = \int_{|z|}^{\sqrt{r^2+z^2}} R^{\alpha+1} dR = \begin{cases}
        \frac{1}{\alpha+2}\left(\left(r^2 + z^2\right)^{\frac{\alpha+2}{2}} - |z|^{\alpha+2}\right), &\alpha \neq -2, \\
        \frac{1}{2}\log\left(\frac{r^2}{z^2} + 1\right), &\alpha = -2.
    \end{cases}
\end{equation}

We can directly reduce Eq. (\ref{eqn:app_I2}) to a combination of simpler integrals:
\begin{equation}
    I_{\alpha,2}(r,z) = \int_0^r \left(\rho^2 + z^2\right)^{\frac{\alpha+2}{2}} - z^2\left(\rho^2 + z^2\right)^{\frac{\alpha}{2}} \, d\rho = J_{\alpha+2}(r,z) - z^2J_{\alpha}(r,z),
\end{equation}
where:
\begin{equation}
    J_\alpha(r,z) = \int_0^r \left(\rho^2+z^2\right)^{\frac{\alpha}{2}} d\rho.
\end{equation}
Integrating by parts, we obtain:
\begin{equation}
    J_\alpha(r,z) = r\left(r^2+z^2\right)^{\frac{\alpha}{2}} - \alpha\int_0^r \left(\rho^2+z^2\right)^{\frac{\alpha-2}{2}}\rho^2 d\rho = r\left(r^2+z^2\right)^{\frac{\alpha}{2}} - \alpha \left(J_\alpha(r,z) - z^2J_{\alpha-2}(r,z)\right),
\end{equation}
resulting in two recursion relationships:
\begin{align}
    J_\alpha(r,z) & = \frac{1}{\alpha+1} \left(r\left(r^2+z^2\right)^{\frac{\alpha}{2}} + \alpha z^2 J_{\alpha-2}(r,z)\right), &\alpha \neq -1, \\
    J_{\alpha}(r,z) & = -\frac{1}{(\alpha+2) z^2} \left(r\left(r^2+z^2\right)^{\frac{\alpha+2}{2}} - (\alpha+3) J_{\alpha+2}(r,z)\right), &\alpha \neq -2.
\end{align}
We can then write the exact expressions using the following base cases:
\begin{align}
    J_{-1}(r,z) & = \int_0^r \frac{1}{\sqrt{\rho^2+z^2}} d\rho = \mbox{arsinh}\left(\frac{r}{z}\right), \\
    J_{-2}(r,z) & = \int_0^r \frac{1}{\rho^2+z^2} d\rho = \frac{1}{z}\arctan\left(\frac{r}{z}\right).
\end{align}

\bibliographystyle{elsarticle-num} 
\bibliography{ref}

\end{document}